\documentclass{article}

\usepackage{amsmath}
\usepackage{amsthm}
\usepackage{amsfonts}
\usepackage{mathrsfs}
\usepackage{amssymb}
\usepackage{bbm}
\usepackage{yhmath}

\usepackage{tikz}
\usetikzlibrary{cd,matrix,arrows,positioning,decorations.pathmorphing}

\usepackage{enumerate}

\usepackage{hyperref}
\hypersetup{colorlinks=true, linkcolor=black, filecolor=black, urlcolor=black, citecolor=black}

\theoremstyle{definition}
\newtheorem{theorem}{Theorem}[subsection]

\newtheorem{lemma}[theorem]{Lemma}
\newtheorem{proposition}[theorem]{Proposition}
\newtheorem{corollary}[theorem]{Corollary}
\newtheorem{remark}[theorem]{Remark}

\def\calG{\mathcal{G}}

\def\calL{\mathcal{L}}
\def\calM{\mathcal{M}}

\def\calO{\mathcal{O}}

\def\calZ{\mathcal{Z}}

\def\frb{\mathfrak{b}}
\def\frc{\mathfrak{c}}

\def\frg{\mathfrak{g}}
\def\frh{\mathfrak{h}}

\def\ffrm{\mathfrak{M}}

\def\bbF{\mathbb{F}}

\def\bbP{\mathbb{P}}
\def\bbQ{\mathbb{Q}}

\def\bbZ{\mathbb{Z}}

\DeclareFontFamily{U}{wncy}{}
\DeclareFontShape{U}{wncy}{m}{n}{<->wncyr10}{}
\DeclareSymbolFont{mcy}{U}{wncy}{m}{n}
\DeclareMathSymbol{\Sha}{\mathord}{mcy}{"58} 

\def\vphi{\varphi}

\def\bs{\backslash}

\DeclareMathOperator{\ox}{\otimes}

\def\comp{\circ}
\def\inj{\hookrightarrow}
\def\surj{\twoheadrightarrow}
\def\isom{\cong}
\def\aisom{\xrightarrow{\sim}}

\DeclareMathOperator{\im}{im}

\DeclareMathOperator{\coker}{coker}

\DeclareMathOperator{\id}{id}

\newcommand{\dlim}{\varinjlim}
\newcommand{\ilim}{\varprojlim}

\DeclareMathOperator{\Fil}{Fil}

\DeclareMathOperator{\Hom}{Hom}

\DeclareMathOperator{\Ext}{Ext}

\DeclareMathOperator{\Gal}{Gal}

\DeclareMathOperator{\Frob}{Frob}

\DeclareMathOperator{\GL}{GL}

\DeclareMathOperator{\End}{End}

\DeclareMathOperator{\Lie}{Lie}

\DeclareMathOperator{\ad}{ad}

\DeclareMathOperator{\gl}{\mathfrak{gl}}

\def\Mod{\text{Mod}}

\def\bb1{\mathbbm{1}}

\DeclareMathOperator{\ur}{ur}

\DeclareMathOperator{\dR}{dR}
\DeclareMathOperator{\cris}{cris}
\DeclareMathOperator{\st}{st}

\DeclareMathOperator{\Spf}{Spf}

\DeclareMathOperator{\Spa}{Spa}

\def\LT{\text{LT}}

\def\et{{\text{\'et}}}

\def\an{\text{an}}
\def\lan{\text{la}}
\def\sm{\text{sm}}

\def\alg{\text{alg}}

\def\hat{\widehat}
\def\bar{\overline}
\def\tilde{\widetilde}

\def\-{\text{-}}

\def\pst{\text{pst}}
\def\WD{\text{WD}}

\def\pro{\text{pro}}

\def\-{\text{-}}

\def\fl{\mathscr{F}\ell}

\def\Sym{\text{Sym}}

\def\GM{\text{GM}}

\def\Nilp{\text{Nilp}}

\def\Dr{\text{Dr}}

\def\geo{\text{geo}}

\def\JL{\text{JL}}

\def\HK{\text{HK}}
\def\DF{\text{DF}}
\def\LL{\text{LL}}

\def\ul{\underline}

\def\ov{\stackrel}

\newcommand{\cf}{\emph{cf.}\ }

\usepackage{geometry}
\usepackage{indentfirst}  

\usepackage{yhmath}

\begin{document}
\title{On the locally analytic $\text{Ext}^1$-conjecture in the $\text{GL}_2(L)$ case}
\date{}
\author{Benchao Su}
\maketitle

\subsubsection*{Abstract}
Let $L$ be a finite extension of $\mathbb{Q}_p$. We calculate the dimension of $\text{Ext}^1$-groups of certain locally analytic representations of $\text{GL}_2(L)$ defined using coherent cohomology of Drinfeld towers of dimension $1$ over $L$. Furthermore, let $\rho_p$ be a $2$-dimensional continuous representation of $\text{Gal}(\bar L/L)$, which is de Rham with parallel Hodge-Tate weights $0,1$ and whose underlying Weil-Deligne representation is irreducible. We prove Breuil's locally analytic $\text{Ext}^1$ conjecture for such $\rho_p$. As an application, we show that the isomorphism class of the multiplicity space $\Pi^{\text{an}}_{\text{geo}}(\rho_p)$ of $\rho_p$ in the pro-\'etale cohomology of Drinfeld towers of dimension $1$ over $L$ considered in \cite{CDN20} uniquely determines the isomorphism class of $\rho_p$.

\tableofcontents

\section{Introduction}
Let $L/\bbQ_p$ be a finite extension, and let $E$ be a sufficiently large finite extension of $\bbQ_p$. The $p$-adic local Langlands program for the group $\GL_2(L)$ aims to construct a correspondence between certain $2$-dimensional $p$-adic representations of $\Gal(\bar L/L)$ over $E$ and certain $p$-adic representations of $\GL_2(L)$ over $E$. In the case $L=\bbQ_p$, such a correspondence is well-established (\cf \cite{breuil2010emerging, colmez2010representations, emerton2011local}). However, when $L\neq\bbQ_p$, the $p$-adic local Langlands program for $\GL_2(L)$ is still very mysterious.

We will mainly focus on the locally analytic aspect. When $L=\bbQ_p$ and the Galois representation is de Rham with regular Hodge-Tate weights, whose underlying Weil-Deligne representation is irreducible, it is known that such a correspondence can be geometrically realized in the cohomology of Drinfeld tower of dimension $1$ over $L$ (\cf \cite{DLB17,CDN20,vanhaecke2024cohomologiesystemeslocauxpadiques}). One of the key ingredients in the proof is the Breuil-Strauch conjecture \cite{DLB17}, which crucially uses the known $p$-adic local Langlands correspondence for $\GL_2(\bbQ_p)$.

When $L\neq\bbQ_p$, the work \cite{CDN20} also provides some structural results on the pro-\'etale cohomology of the Drinfeld tower of dimension $1$ over $L$. Let us introduce some notation. Given a $2$-dimensional continuous $E$-linear representation $V$ of $\Gal(\bar L/L)$, which is de Rham with parallel Hodge-Tate weight $0,1$ (our convention is that the Hodge-Tate weight of the $p$-adic cyclotomic character is $1$), we can associate the following objects: 
\begin{itemize}
    \item Let $M$ be the Deligne-Fontaine module over $E\ox_{\bbQ_p}L_0^{\ur}$ associated to $V$ via Fontaine's recipe. We assume $D_{\pst}(V)\isom M[-1]$. Here $M[-1]$ is a twist of $M$ such that $\vphi_{M[-1]}=p^{-1}\vphi_M$.
    \item Let $\WD(M)$ be the Weil-Deligne representation associated to $M$. This is a $2$-dimensional $E$-linear representation of $\WD_L$, the Weil-Deligne group of $L$. In the most part of the article we will assume that $\WD(M )$ is irreducible.
    \item Let $\text{LL}(M )$ be the smooth irreducible representation of $\GL_2(L)$ corresponding to $\WD(M )$ via the local Langlands correspondence. If $\WD(M )$ is irreducible, the smooth $\GL_2(L)$-representation $\LL(M )$ is supercuspidal.
    \item Let $\JL(M )$ be the smooth irreducible representation of $D_L^\times$ corresponding to $\LL(M )$ via the Jacquet-Langlands correspondence. Here $D_L$ is a non-split quaternion algebra over $L$. If $\LL(M )$ is supercuspidal, the action of $D_L^\times$ on $\JL(M )$ does not factor through the reduced norm map. 
\end{itemize}

The key geometric objects are the Drinfeld tower of dimension $1$. Let $\Omega:=\bbP^1_C\bs \bbP^1(L)$ be the Drinfeld's upper half-plane, viewed as an adic space over $C$, where $C$ is the completion of an algebraic closure of $L$. This space $\Omega$ has a sequence of \'etale Galois coverings $\{\calM_{\Dr,n}\}_{n\ge 0}$, with Galois group $D_L^\times/(1+\varpi^n\calO_{D_L})$. Here, we fix a uniformizer $\varpi$ of $L$, and we view $\varpi$ as an element of $D_L^\times$ by identifying the center of $D_L^\times$ with $L^\times$. We denote by $\calO_{D_L}$ the maximal order of $D_L$. Also each space $\calM_{\Dr,n}$ carries an action of $\GL_2(L)$, commuting with the action of $D_L^\times$. From now on we assume $\WD(M)$ is irreducible (so that we can find some related representations on the Drinfeld tower of dimension $1$). Define 
\begin{align}
    H^1_{\pro\et}[M]:=\Hom_{E[D_L^\times]}(\JL(M),E\ox_{\bbQ_p}\dlim_n H^1_{\pro\et}(\calM_{\Dr,n},\bbQ_p(1))),
\end{align}
where $\bbQ_p(1)$ is the Tate twist. This is an ind-Fr\'echet space over $E$, carrying an action of $\GL_2(L)\times \Gal(\bar L/L)$. In order to describe it, we also need to define some objects from the de Rham complex of $\calM_{\Dr,n}$ for $n\ge 0$: 
\begin{align}
    \calO[M]&:=\Hom_{E[D_L^\times]}(\JL(M),E\ox_{\bbQ_p}\dlim_n H^0(\calM_{\Dr,n},\calO_{\calM_{\Dr,n}}))^{\Gal(\bar L/L)},\\
    \Omega^1[M]&:=\Hom_{E[D_L^\times]}(\JL(M),E\ox_{\bbQ_p}\dlim_n H^0(\calM_{\Dr,n},\Omega^1_{\calM_{\Dr,n}}))^{\Gal(\bar L/L)}.
\end{align}
They are all ind-Fr\'echet spaces over $E\ox_{\bbQ_p}L$, whose strong dual are locally $\bbQ_p$-analytic representations of $\GL_2(L)$. One can show that they are all non-zero (Corollary \ref{cor:nonzero}). The differential map on $\calM_{\Dr,n}$ gives an exact sequence
\begin{align}
    0\to \calO[M]\to \Omega^1[M]\to M_{\dR}\ox_{E}\LL(M)'\to 0,
\end{align}
where $\LL(M)'$ is the strong dual of $\LL(M)$ (equipped with ind-finite dimensional topology), and $M_{\dR}:=(C\ox_{L_0^{\ur}}M)^{\Gal(\bar L/L)}$ is a rank $2$ free module over $L\ox_{\bbQ_p}E$, viewed as an $E$-algebra via $E\to L\ox_{\bbQ_p}E,x\mapsto 1\ox x$. Using the decomposition $L\ox_{\bbQ_p}E\isom \prod_{\sigma\in \Sigma}E$, we can define 
\begin{align}\label{5}
    0\to \calO_\sigma[M]\to \Omega^1_\sigma[M]\to M_{\dR,\sigma}\ox_{E}\LL(M)'\to 0,
\end{align}
where the subscript $(-)_\sigma$ denotes the functor $-\ox_{L\ox_{\bbQ_p}E,\sigma\ox \id}E$. The terms in (\ref{5}) are naturally modules over the locally $\sigma$-analytic distribution algebra $D_\sigma(\GL_2(L),E)$, which is the strong dual of the space of locally $\sigma$-analytic functions on $\GL_2(L)$. 

Up to now, all these $\GL_2(L)$-representations only depend on the Deligne-Fontaine module $M$. Here is how the Hodge filtration of $V$ enters in. Recall that to the $p$-adic Galois representation $V$ we can associate a filtered $L\ox_{\bbQ_p}E$-vector space $D_{\dR}(V)$, whose filtration is called the Hodge filtration of $V$. The decomposition $L\ox_{\bbQ_p}E\isom \prod_{\sigma\in \Sigma}E$ also gives a decomposition $D_{\dR}(V)\isom \bigoplus_{\sigma\in\Sigma}D_{\dR,\sigma}(V)$. From the Hyodo-Kato isomorphism (applied to $M_{\dR,\sigma}$) and general theory of $p$-adic Galois representations (applied to $D_{\dR,\sigma}(V)$), we know that there is an isomorphism of $E$-vector spaces $M_{\dR,\sigma}\isom D_{\dR,\sigma}(V)$ by using their underlying structure of Weil-Deligne representations. As $\WD(M)$ is irreducible, we know that such an isomorphism is unique up to a scalar. For each $\sigma\in \Sigma$, define $\calL_\sigma=\Fil^0D_{\dR,\sigma}(V)$. This is a $1$-dimensional subspace inside $D_{\dR,\sigma}(V)$. From the isomorphism $M_{\dR,\sigma}\isom D_{\dR,\sigma}(V)$, we also view $\calL_\sigma$ as a subspace of $M_{\dR,\sigma}$. Let $W_{M,\calL_\sigma}$ be the inverse image of $\calL_\sigma\ox_E\LL(M)'\subset M_{\dR,\sigma}\ox_E\LL(M)'$ inside $\Omega^1_\sigma[M]$, and we denote $\Pi^{\an}_{\geo,\sigma}(V)$ by the strong dual of $W_{M,\calL_\sigma}$. This is a locally $\sigma$-analytic representation of $\GL_2(L)$ over $E$. Besides, define $\Pi^{\an}_{\geo}(V)$ to be the strong dual of $\Hom_{E[\Gal(\bar L/L)]}(V,H^1_{\pro\et}[M])$. By \cite[Proposition 2.10]{CDN20},  there is a $\GL_2(L)$-equivariant isomorphism
\begin{align}
    \Pi^{\an}_{\geo}(V)\isom \bigoplus_{\sigma\in \Sigma,\LL(M)}\Pi^{\an}_{\geo,\sigma}(V)
\end{align}
where the right hand side is an amalgamated sum of $W_{M,\calL_\sigma}$ over $\LL(M)$ for varying $\sigma$. So one may imagine the structure of $\Pi^{\an}_{\geo}(V)$ as follows: 
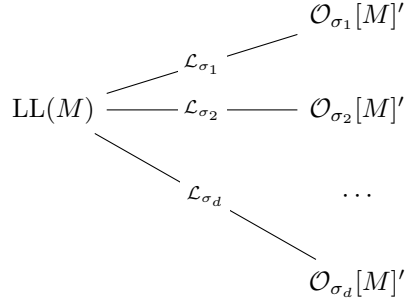
\[\begin{tikzcd}
	&&& {\mathcal{O}_{\sigma_1}[M]'} \\
	{\text{LL}(M)} &&& {\mathcal{O}_{\sigma_2}[M]'} \\
	&&& \cdots \\
	&&& {\mathcal{O}_{\sigma_d}[M]'}
	\arrow["{\mathcal{L}_{\sigma_1}}"{description}, no head, from=2-1, to=1-4]
	\arrow["{\mathcal{L}_{\sigma_2}}"{description}, no head, from=2-1, to=2-4]
	\arrow["{\mathcal{L}_{\sigma_d}}"{description}, no head, from=2-1, to=4-4]
\end{tikzcd}\]
where we lable embeddings in $\Sigma$ by $\sigma_1,...,\sigma_d$ with $d=[L:\bbQ_p]$, and $\mathcal{O}_{\sigma_i}[M]'$ denotes the strong dual of $\mathcal{O}_{\sigma_i}[M]$. Fix an embedding $\sigma\in \Sigma$. Here are some natural questions:
\begin{itemize}
    \item Does $W_{M,\calL_\sigma}\isom W_{M,\calL_\sigma'}$ (as $\GL_2(L)$-representations) imply $\calL_\sigma=\calL_\sigma'$?
    \item Given any non-split extension $W$ of $\LL(M)'$ by $\calO_\sigma[M]$, does there exist some $\calL_\sigma$ such that $W\isom W_{M,\calL_\sigma}$? In other words, does there exist an embedding $W\inj \Omega^1_\sigma[M]$? This amounts to asking whether $\Omega^1_\sigma[M]$ is the universal extension of $\LL(M)'$ by $\calO_\sigma[M]$.
\end{itemize}
In fact, these questions are closely related to Breuil's locally analytic $\Ext^1$-conjecture in the $\GL_2(L)$-case (\cf\cite{Bre19}). When $L=\bbQ_p$, Ding proves this conjecture \cite{Din22}, which is based on the $p$-adic local Langlands correspondence for $\GL_2(\bbQ_p)$ and computations of $(\vphi,\Gamma)$-modules. In this paper, we answer these questions and give a different proof which works for arbitrary base field $L$. 
\begin{theorem}\label{thm:intromain}
\begin{enumerate}[(i)]
    \item (Corollary \ref{cor:Ext1onlyonDr}) We have $\dim_E\Ext^1_{D_\sigma(\GL_2(L),E)}(\LL(M)',\calO_\sigma[M])=2$.
    \item (Theorem \ref{thm:inj}) The $D_\sigma(\GL_2(L),E)$-module $\Omega^1_\sigma[M]$ is the universal extension of $\LL(M)'$ by $\calO_\sigma[M]$.
    \item (Theorem \ref{thm:Ext1Dr}) There exists an isomorphism of $E$-vector spaces 
    \begin{align}
        \WD(M[-1])\aisom \Ext^1_{D_\sigma(\GL_2(L),E)}(\LL(M)',\calO_\sigma[M]),
    \end{align}
    which only depends on $M$, such that for any $2$-dimensional $E$-linear de Rham representation $V$ of parallel Hodge-Tate weights $0,1$, and $D_{\pst}(V)\isom M[-1]$, under the natural isomorphism (unique up to a scalar)
    \begin{align}
        D_{\dR,\sigma}(V)\isom D_{\pst}(V)_\sigma\aisom \WD(M[-1])\aisom \Ext^1_{D_\sigma(\GL_2(L),E)}(\LL(M)',\calO_\sigma[M]),
    \end{align}
    the $E$-line $\Fil^0D_{\dR,\sigma}(V)$ is mapped to the $E$-line generated by the extension class given by $W_{M,\Fil^0D_{\dR,\sigma}(V)}$.
\end{enumerate}
\end{theorem}
For normalizations of these twists, we refer the reader to Section \ref{sec:prelim}. Actually, we prove the above results for general regular algebraic weights, where one replaces the representations defined using the de Rham complex by general dual BGG complexes.

\begin{remark}
Using a similar method one can also obtain the following result by reversing the target and the source:
\begin{align}
    \dim_E \Ext^1_{D_\sigma(\GL_2(L),E)}(\calO_\sigma[M],\LL(M)')=2.
\end{align}
\end{remark}


A key ingredient in the proof is to compute the wall-crossing of the $D_\sigma(\GL_2(L),E)$-module $\calO_\sigma[M]$. Roughly speaking, wall-crossing functors for locally analytic representations \cite{JLS22, Din24} are natural generalizations of wall-crossing functors for Lie algebra representations, for example \cite[Chapter 7]{Hum94}. The wall-crossing functors are helpful for analyzing the structure of some locally analytic representations and are also useful for producing new ones, see for example \cite[Theorem 1.3.7]{breuil2024splittingexplicitingrhamcomplex}. As $\calM_{\Dr,n}$ for $n\ge 0$ are \'etale over the period domain $\bbP^1$, and the wall-crossing functors commute with taking cohomology, essentially one just needs to compute the wall-crossing of the equivariant line bundles on $\bbP^1$, which we handle it by explicit calculations. It turns out that the wall-crossing of $\calO_\sigma[M]$ has some nice properties, basically coming from the following key points:
\begin{enumerate}[(i)]
    \item functoriality of the wall-crossing functor,
    \item adjointness of translation functors (as wall-crossing functors are compositions of translation functors),
    \item the fact that the wall-crossing of $\calO_\sigma[M]$ follows the same pattern as the wall-crossing of the structure sheaf of $\bbP^1$ (also the wall-crossing of the irreducible dual Verma module in the principal block of the BGG category $\calO$).
\end{enumerate}
Therefore, once we relate $\Omega^1_\sigma[M]$ with the wall-crossing of $\calO_\sigma[M]$ (Proposition \ref{prop:injsurj}), we can show that $\Omega^1_\sigma[M]$ is the universal extension of $\LL(M)'$ by $\calO_\sigma[M]$.

Finally, after additional steps on describing endomorphisms of $\calO_\sigma[M]$ (Lemma \ref{lem:wcsemisimple}), we prove that the pro-\'etale cohomology of Drinfeld tower of dimension $1$ over $L$ realizes a one-to-one correspondence between the following two types of objects:
\begin{itemize}
    \item $2$-dimensional de Rham representation of $\Gal(\bar L/L)$ over $E$, of parallel Hodge-Tate weights $0,1$, whose underlying Weil-Deligne representation is irreducible, up to $E$-linear isomorphism;
    \item certain locally $\bbQ_p$-analytic representations of $\GL_2(L)$ over $E$ that are constructed using the de Rham complex of Drinfeld tower of dimension $1$ over $L$, up to continuous $E$-linear isomorphism. 
\end{itemize}
We note that this result is also predicted by the conjecture (GAL) in \cite[2.3]{Bre19}.
\begin{theorem}[Corollary \ref{cor:MainDr}]
Let $V$ be as in Theorem \ref{thm:intromain}. The $E$-linear isomorphism class of $\Pi^{\an}_{\geo}(V)$ determines the Hodge filtration of $V$.
\end{theorem}

\subsubsection*{Acknowledgements}
We thank Yiwen Ding, Zicheng Qian, James Taylor,  Arnaud Vanhaecke, and Liang Xiao for helpful discussions. We thank Arnaud Vanhaecke for his thorough reading of this article, which significantly improved the content. Finally, we sincerely thank our Ph.D. advisor Yiwen Ding for his patience, for pointing out some crucial gaps in early versions of the draft, and for providing many key comments, from which we benefit greatly.

\subsubsection*{Notations and conventions}
Let $L$ be a finite extension of $\bbQ_p$ (this is the base field), and let $E$ be a finite extension of $\bbQ_p$ (this is the coefficient field) such that $[L:\bbQ_p]=|\Hom_{\bbQ_p}(L,E)|$. We will always assume $E$ to be sufficiently large. Let $L_0\subset L$ be the maximal unramified extension of $\bbQ_p$. Let $L_0^{\ur}$ be the maximal unramified subextension of $L_0$, and let $\breve{L}_0$ be the $p$-adic completion of $L_0^{\ur}$. We always view $L_0$ as a subfield of $L$ via the natural embedding. Fix an algebraic closure $\bar L$ of $L$.  Let $\text{W}_L$ be the Weil group of $L$, and let $\WD_L$ be the Weil-Deligne group of $L$. Let $C:=\hat{\bar L}$ be the $p$-adic completion of $\bar L$. Let $B_{\dR}$ be the field of $p$-adic periods. Let $B_{\cris}\subset B_{\st}\subset B_{\dR}$ be the rings of crystalline (resp. semi-stable) periods, where we fix a uniformizer $\varpi$ of $L$, and we fix a branch of the $p$-adic logarithm $\log_\varpi$ such that $\log_\varpi(\varpi)=0$. Our convention is that the Hodge-Tate weight of the $p$-adic cyclotomic character of $\Gal(\bar \bbQ_p/\bbQ_p)$ is $1$.

For $\pi_1,\cdots,\pi_n$, a list of locally $\bbQ_p$-analytic representations of $\GL_2(L)$, we write $\pi_1-\pi_2-\cdots-\pi_n$ to denote a successive extension of $\pi_1,\cdots,\pi_n$ such that for each $k$, $\pi_k-\pi_{k+1}-\cdots-\pi_n$ is a non-split extension of $\pi_{k+1}-\cdots-\pi_n$ by $\pi_k$.

\section{Wall-crossing on Drinfeld towers of dimension $1$}
In this section, we will recall some equivariant sheaves on Drinfeld towers of dimension $1$ over $L$ , and apply the wall-crossing method to these sheaves and the cohomology of these sheaves. Finally, we determine the dimension of some extension groups.

\subsection{Drinfeld towers of dimension $1$}
Let $L$ be a finite extension of $\bbQ_p$, with $\calO_L\subset L$ the ring of integers in $L$, and residue field $\bbF_q$. We fix a uniformizer $\varpi$ of $L$. Let $\breve{L}$ be the completion of the maximal unramified extension of $L$, and $\calO_{\breve{L}}$ be the ring of integers of $\breve{L}$. Let $C$ be the completion of an algebraic closure of $L$. Let $D_L$ be a non-split quaternion algebra over $L$, and $\calO_{D_L}$ the maximal order of $D_L$. First we recall the definition of special formal $\calO_{D_L}$-modules. Let $S$ be a scheme over $\calO_{\breve{L}}$. A special formal $\calO_{D_L}$-module over $S$ is a pair $(X,\sigma)$ consisting of a formal group $X$ over $S$ and a ring homomorphism $\sigma:\calO_{D_L}\inj \End(X)$ such that the action of $\calO_L$ on $\Lie(X)$ induced by $\sigma$ coincides with that induced by $\calO_L\inj\calO_{\breve{L}}\to\calO_S$, and $\Lie(X)$ is a free $\calO_{L_2}\ox_{\calO_L}\calO_S$-module of rank $1$, where $\calO_{L_2}\subset \calO_{D_L}$ is the integral ring of the quadratic unramified extension inside $D_L$.

Let $H_1$ be a special formal $\calO_{D_L}$-module over $\bar\bbF_q$ of $\calO_L$-height $4$, which is unique up to quasi-isogeny. Here $q$ is the cardinality of the residue field of $L$. The group of $\calO_{D_L}$-quasi-isogenies of $H_1$ is $\text{QIsog}_{\calO_{D_L}}(H_1)\isom \GL_2(L)$. Let $\Nilp_{\calO_{\breve{L}}}$ be the category of $\calO_{\breve{L}}$-algebras on which $\varpi$ is nilpotent. Let $\ffrm_{\Dr,0}$ be the functor from $\Nilp_{\calO_{\breve{L}}}$ to the category of sets, which assigns an $\calO_{\breve{L}}$-algebra $R$ to the set of equivalence classes of pairs $(G,\rho)$, where 
\begin{itemize}
    \item $G$ is a $p$-divisible special formal $\calO_{D_L}$-module over $R$.
    \item $\rho:H_1\ox_{\bar\bbF_q}R/\varpi\dasharrow G\ox_R R/\varpi$ is an $\calO_{D_L}$-equivariant quasi-isogeny.
\end{itemize}
Two pairs $(G,\rho)$ and $(G',\rho')$ are equivalent if the quasi-isogeny $\beta'\comp\beta^{-1}:G\ox_R R/\varpi\dasharrow G'\ox_R R/\varpi$ lifts to an isomorphism $\tilde{\beta'\comp\beta^{-1}}:G\aisom G'$ over $R$. Drinfeld \cite{Dr76} showed that the functor $\ffrm_{\Dr,0}$ is represented by a formal scheme over $\Spf \calO_{\breve{L}}$. It is equipped with a $\GL_2(L)$-action: for $g\in \GL_2(L)\isom \text{QIsog}_{\calO_{D_L}}(H_1)$, it acts as $(G,\rho)\mapsto (G,\rho\comp g^{-1})$. Let $\breve\calM_{\Dr,0}$ be the adic generic fiber of $\ffrm_{\Dr,0}$. The space $\breve{\calM}_{\Dr,0}$ is equipped with a canonical Weil descent datum relative to $\breve L/L$ and we have an isomorphism of adic spaces:
\begin{align}
    \breve{\calM}_{\Dr,0}\aisom(\Omega_L\times_{\Spa(L,\calO_L)}\Spa(\breve{L},\calO_{\breve{L}}))\times\bbZ
\end{align}
compatible with the Weil-descent datum, see \cite[Theorem 3.72]{RZ96}. Here $\Omega_L=\bbP^{1,\ad}_L\bs\bbP^1(L)\subset\bbP^1_L$ is the Drinfeld upper half plane over $L$ and $\bbP^{1,\ad}_L$ is the adic space associated to the projective line $\bbP^1_L$ over $L$. This isomorphism is given by the period map which we will define later.  By considering the Hodge filtration of the Dieudonn\'e module of the universal $p$-divisible group on $\breve{\calM}_{\Dr,0}$, we get the Gross-Hopkins period map $\breve\pi_{\GM}$ on $\breve{\calM}_{\Dr,0}$, which is given by
\begin{align}
    \breve{\calM}_{\Dr,0}\to\Omega_L\times_{\Spa(L,\calO_L)}\Spa(\breve{L},\calO_{\breve{L}})\to \bbP^{1,\text{ad}}_L\times_{\Spa(L,\calO_L)}\Spa(\breve{L},\calO_{\breve{L}})
\end{align}
where the first map is obtained from the isomorphism by projecting to the first factor, and the second map follows by the natural inclusion $\Omega_L\subset\bbP^{1,\text{ad}}_L$.

There are some \'etale coverings of $\breve{\calM}_{\Dr,0}$ obtained by trivializing the universal $p$-adic Tate module on $\breve{\calM_{\Dr,0}}$. Let $(\calG',\rho')$ be the universal $p$-divisible special formal $\calO_{D_L}$-module on $\breve\calM_{\Dr,0}$. The $p$-adic Tate module of $\calG'$ defines a $\bbZ_p$-local system equipped with an action of $\calO_{D_L}$ on the \'etale site of $\breve\calM_{\Dr,0}$, which we denote by $V_{\Dr}$. For $n\ge 0$, we get an \'etale Galois covering $\breve\calM_{\Dr,n}$ of $\breve\calM_{\Dr,0}$ with Galois group $(\calO_{D_L}/\varpi^n\calO_{D_L})^\times=\calO_{D_L}^\times/(1+\varpi^n\calO_{D_L})$, by considering $\calO_{D_L}$-equivariant trivializations of $V_{\Dr}/\varpi^n V_{\Dr}$. These spaces are equipped with $\GL_2(L)\times D_L^\times$-actions. For $m\ge n\ge 0$, the covering maps $\breve\calM_{\Dr,m}\to\breve\calM_{\Dr,n}$ are $\GL_2(L)\times D_L^\times$-equivariant. Moreover, for $n\ge 0$, the spaces $\breve{\calM}_{\Dr,n}$ are equipped with canonical Weil descent data relative to $\breve L/L$, which are not effective but become effective by quotient ${\varpi^{\bbZ}}$, with $\varpi$ viewed as a central element in $\GL_2(L)$. Fix a $\bbQ_p$-embedding $\sigma:L\inj E$ and let $\Sigma_n$ be the adic space over $E$ obtained by base-changing the $L$-model of $\breve\calM_{\Dr,n}/{\varpi^{\bbZ}}$ via $\sigma$. The map $\breve{\pi}_{\GM}$ induces a natural map $\pi_{\GM}:\Sigma_0\to \fl$ with $\fl$ the adic space associated to $\bbP^1_L\times_{L,\sigma} E$. For any $n\ge 0$, let $\pi_n:\Sigma_n\to \Sigma_0$ be the projection map and denote the composition map $\pi_{\GM}\comp \pi_n:\Sigma_n\to \fl$ by $\pi_{\GM,n}$. The map $\pi_{\GM,n}$ is compatible with the $\GL_2(L)\times D_L^\times$-action, with the natural action of $\GL_2(L)$ on $\fl$ (via $\sigma$) and trivial action of $D_L^\times$ on $\fl$.

Using the period map $\pi_{\GM,n}$, we can pull back some equivariant sheaves on $\fl$ to $\Sigma_n$. Let $\bar B\subset \GL_2$ be the Borel subgroup consisting of lower triangular matrices, and we identify $\fl$ with the base change to $E$ via $\sigma$ of the adic space associated with $\bar B\bs \GL_2$. For $(a,b)\in\bbZ^2$, let $\omega_{\fl}^{(a,b)}$ be the $\GL_2(L)$-equivariant line bundle on $\fl$ associated to the $\bar B$-representation $\left(\begin{matrix}
    x&0\\y&z
\end{matrix}  \right)\mapsto \sigma(x)^a\sigma(z)^b$. Under this normalization, one directly checks that $H^0(\fl,\omega_{\fl}^{(1,0)})$ is the $\sigma$-standard representation of $\GL_2(L)$, and we have a canonical $\GL_2(L)$-equivariant isomorphism $\Omega_{\fl}^1\isom \omega_{\fl}^{(-1,1)}$. Here, by the $\sigma$-standard representation of $\GL_2(L)$ we mean the base change to $E$ via $\sigma$ of the standard action of $\GL_2(L)$ on a $2$-dimensional $L$-vector space. Define $\omega_{\Sigma_n}^{(a,b)}:=\pi_{\GM,n}^*\omega_{\fl}^{(a,b)}$ with $\calO_{\Sigma_n}=\omega_{\Sigma_n}^{(0,0)}$ the structure sheaf of $\Sigma_n$ and $\Omega^1_{\Sigma_n}=\omega_{\Sigma_n}^{(-1,1)}$ the sheaf of differentials on $\Sigma_n$. As $\calO_{\Sigma_n}$-modules, these sheaves $\omega_{\Sigma_n}^{(a,b)}$ for $(a,b)\in\bbZ^2$ are all isomorphic to the structure sheaf $\calO_{\Sigma_n}$, but the isomorphisms are not compatible with the $\GL_2(L)$-action.

Finally we introduce some notation on Lie algebras. Fix a $\bbQ_p$-embedding $\sigma:L\inj E$. Let $\frg_\sigma:=\gl_2(L)\ox_{L,\sigma}E$, and $\calZ(\frg_\sigma)$ be the center of $U(\frg_{\sigma})$.  The center $\calZ(\frg_\sigma)$ is generated as an $E$-algebra by $z=\left( \begin{matrix} 1&0\\0&1\end{matrix} \right)\in \frg_\sigma$, and the Casimir operator $\frc=u^+u^-+u^-u^++\frac{1}{2}h^2\in U(\frg_\sigma)$, where $u^+=\left( \begin{matrix} 0&1\\0&0\end{matrix} \right)\in\frg_\sigma$, $u^-=\left( \begin{matrix} 0&0\\1&0\end{matrix} \right)\in\frg_\sigma$, and $h=\left( \begin{matrix} 1&0\\0&-1\end{matrix} \right)\in\frg_\sigma$. The action of $\calZ(\frg_\sigma)$ on $\omega_{\Sigma_n}^{(a,b)}$ is semisimple, for $(a,b)\in\bbZ^2$. Let $\frb_{\sigma}$ be the Borel subalgebra consisting of upper triangular matrices, and let $\frh_{\sigma}$ be the Cartan subalgebra consisting of diagonal matrices. For an integral weight $\mu$ of $\frh_\sigma$, we denote by $\tilde{\chi}_\mu$ the infinitesimal character of $\calZ(\frg_\sigma)$ acting on $U(\frg_{\sigma})\ox_{U(\frb_\sigma)}\mu$. Then $\omega_{\Sigma_n}^{(a,b)}$ has infinitesimal character $\tilde{\chi}_{(a,b)}$. Indeed, $\pi_{\GM,n}$ is an \'etale map and $\calZ(\frg_\sigma)$ acts on $\omega_{\fl}^{(a,b)}$ via $\tilde{\chi}_{(a,b)}$.

\subsection{Translations of equivariant sheaves}
In the classical case, translation functors have been widely applied in the study of representations of semisimple Lie algebras. There exists recent work \cite{JLS22,Din24} applying translation functors to study locally analytic representations. In this section, we adapt these translation functors to equivariant sheaves on Drinfeld towers of dimension $1$ over $L$.

We recall the definition of translation functors for locally $\sigma$-analytic representations of $G=\GL_2(L)$. Let $D_\sigma(G,E)$ be the $E$-valued locally $\sigma$-analytic distribution algebra on $G$, which is the strong dual of the space of $E$-valued locally $\sigma$-analytic functions on $G$. Let $\mu,\lambda\in\Lambda$ be two integral weights of $\frh_\sigma$. Let $\Mod_{\tilde{\chi}_\mu}(D_\sigma(G,E))$ be the category of left $D_\sigma(G,E)$-modules with generalized infinitesimal character $\tilde{\chi}_\mu$.  If $\mu$ is dominant with respect to $\frb_\sigma$, let $V^{\mu}$ be the $\sigma$-algebraic representation of $\GL_2(L)$ over $E$, with highest weight $\mu$ with respect to $\frb_\sigma$. We define the translation functor from weight $\mu$ to weight $\lambda$ as follows:
\begin{align}
    T_\mu^\lambda:\Mod_{\tilde{\chi}_\mu}(D_\sigma(G,E))&\to \Mod_{\tilde{\chi}_\lambda}(D_\sigma(G,E))\\
    M&\mapsto (V^{\bar\nu}\ox_EM)\{\calZ(\frg_\sigma)=\tilde{\chi}_\lambda\}
\end{align}
where $\bar\nu$ is the unique dominant weight in the orbit of the usual action of the Weyl group on $\lambda-\mu$, and $\{-\}$ denotes the generalized eigenspace. By \cite[Lemma 2.4.6, Theorem 2.4.7]{JLS22}, $T_\mu^\lambda$ is an exact functor, with left and right adjoint $T_\lambda^\mu$.

We apply these translation functors to the equivariant sheaves $\omega_{\Sigma_n}^{(a,b)}$ on $\Sigma_n$ we defined before. Let $\mu=(a,b)$ be an integral weight. By \cite[Lemma 2.3.4]{RC1}, for any affinoid open subset $U\subset\Sigma_n$, there exists a sufficiently small open compact subgroup $G_U$ of $\GL_2(L)$ such that the action of $G_U$ on $\omega_{\Sigma_n}^{(a,b)}(U)$ is locally analytic. Moreover, by \cite[Corollary 2.2.6]{RC1}, the action of $G_U$ on $\omega_{\Sigma_n}^{(a,b)}(U)$ is analytic after shrinking $G_U$ if necessary. Besides, as $L$ acts on $\omega_{\Sigma_n}^{(a,b)}$ via $\sigma$, the action of $G_U$ on $\omega_{\Sigma_n}^{(a,b)}(U)$ is locally $\sigma$-analytic. Let $\{U_i\}$ be a Stein covering of $\Sigma_n$. Taking limit on $U_i$, we deduce that $H^0(\Sigma_n,\omega_{\Sigma_n}^{(a,b)})=\ilim_i H^0(U_i,\omega_{\Sigma_n}^{(a,b)})$ is a separately continuous $D_\sigma(G,E)$-module. See also \cite[Remarque 3.4]{DLB17}.  As $\pi_{\GM,n}:\Sigma_n\to \fl$ is an \'etale map and $\calZ(\frg_\sigma)$ acts on $\omega_{\fl}^{(a,b)}$ via $\tilde{\chi}_{(a,b)}$, we see that $\calZ(\frg_\sigma)$ also acts on $\omega_{\Sigma_n}^{(a,b)}$ via $\tilde{\chi}_{(a,b)}$. Let $\lambda$ be an integral weight. As the translation functor $T_\mu^\lambda$ is exact, we deduce that $T_\mu^\lambda(\omega_{\Sigma_n}^{(a,b)}(U))$ for $U$ varying in affinoid open subsets of $\Sigma_n$ defines a sheaf on $\Sigma_n$, by directly checking the sheaf condition. We denote this sheaf by $T_\mu^\lambda\omega_{\Sigma_n}^{(a,b)}$.

We want to show some basic properties of $T_\mu^\lambda\omega_{\Sigma_n}^{(a,b)}$. We start with a lemma.
\begin{lemma}\label{lem:basic}
Let $V$ be a $\sigma$-algebraic representation of $G$ over $E$. For $(a,b)\in\bbZ^2$, the sheaf $V\ox_E\omega_{{\Sigma_n}}^{(a,b)}$ has a filtration with graded pieces given by equivariant locally free sheaves of rank $1$, of the form $\omega_{{\Sigma_n}}^{(a',b')}$ for some $(a',b')\in\bbZ^2$.
\begin{proof}
After twisting by determinants and invertible sheaves, we may assume that $(a,b)=(0,0)$ and $V=V^{(k,0)}\isom \Sym^k V^{(1,0)}$, with $V^{(1,0)}$ the $\sigma$-standard representation of $G$. When $k=1$, the universal flag on $\fl$ gives an exact sequence 
\begin{align}
    0\to \omega_{\fl}^{(0,1)}\to V^{(1,0)}\ox_E\calO_{\fl}\to \omega_{\fl}^{(1,0)}\to 0.
\end{align}
By pulling back to $\Sigma_n$, we get an exact sequence 
\begin{align}
    0\to \omega_{\Sigma_n}^{(0,1)}\to V^{(1,0)}\ox_E\calO_{\Sigma_n}\to \omega_{\Sigma_n}^{(1,0)}\to 0.
\end{align}
For general cases, by taking symmetric powers of the above exact sequence over $\calO_{{\Sigma_n}}$, we deduce that $V\ox_E\calO_{{\Sigma_n}}$ is a  successive extension of $\omega_{{\Sigma_n}}^{(a,k-a)}$ with $a=0,1,...,k$.
\end{proof}
\end{lemma}

\begin{lemma}\label{lem:T0kOfl}
Let $k\ge 0$ be an integer. Then $T_{(0,0)}^{(k,0)}\calO_{{\Sigma_n}}\isom \omega_{{\Sigma_n}}^{(k,0)}$ and $T_{(0,0)}^{(k,0)}\Omega^1_{{\Sigma_n}}\isom \omega_{{\Sigma_n}}^{(-1,k+1)}$.
\begin{proof}
From the proof of Lemma \ref{lem:basic}, we see $V^{(k,0)}\ox\calO_{{\Sigma_n}}$ has a filtration by $\GL_2(L)$-equivariant vector bundles with each graded piece of the form $\omega_{{\Sigma_n}}^{(a,k-a)}$ with $a=0,1,...,k$. Therefore, after we apply $\{\calZ(\frg_\sigma)=\chi_{(k,0)}\}$ to $V^{(k,0)}\ox\calO_{{\Sigma_n}}$, we get $\omega_{{\Sigma_n}}^{(k,0)}$. Similarly one obtains the second claim.
\end{proof}
\end{lemma}

\begin{proposition}\label{prop:basicoftranslation}
Let $n\ge 0$ be an integer, $\lambda,\mu$ be two integral weights. Write $\mu=(a,b)\in\bbZ^2$.
\begin{enumerate}[(i)]
    \item The sheaf $T_\mu^\lambda\omega_{\Sigma_n}^{(a,b)}$ is a coherent sheaf on $\Sigma_n$, and $H^i(\Sigma_n,T_\mu^\lambda\omega_{\Sigma_n}^{(a,b)})=0$ for $i\ge 1$.
    \item We have a $\GL_2(L)\times D_L^\times$-equivariant isomorphism:
    \begin{align}
        H^0(\Sigma_n,T_\mu^\lambda\omega_{\Sigma_n}^{(a,b)})\isom T_\mu^\lambda H^0(\Sigma_n,\omega_{\Sigma_n}^{(a,b)}).
    \end{align}
\end{enumerate}
\begin{proof}
By Lemma \ref{lem:basic}, we deduce that $T_\mu^\lambda\omega_{\Sigma_n}^{(a,b)}$ is a successive extension of line bundles of the form $\omega_{\Sigma_n}^{(a',b')}$ with $(a',b')\in\bbZ^2$. This shows that $T_\mu^\lambda\omega_{\Sigma_n}^{(a,b)}$ is a coherent sheaf on $\Sigma_n$. By Kiehl's vanishing theorem \cite{Kie67}, $H^i(\Sigma_n,T_\mu^\lambda\omega_{\Sigma_n}^{(a,b)})=0$ for $i\ge 1$. Let $\{U_i\}_{i\ge 0}$ be a Stein covering of $\Sigma_n$. As $T^\lambda_\mu$ is biadjoint to $T^\mu_\lambda$, $T^\lambda_\mu$ is exact. From this we deduce that $T_\mu^\lambda$ commutes with taking cohomology:
\begin{align}
    H^0(\Sigma_n,T_\mu^\lambda\omega_{\Sigma_n}^{(a,b)})\isom \ilim_i H^0(U_i,T_\mu^\lambda\omega_{\Sigma_n}^{(a,b)})\isom \ilim_i T_\mu^\lambda H^0(U_i,\omega_{\Sigma_n}^{(a,b)})\isom T_\mu^\lambda \ilim_i H^0(U_i,\omega_{\Sigma_n}^{(a,b)})\isom T_\mu^\lambda H^0(\Sigma_n,\omega_{\Sigma_n}^{(a,b)}).
\end{align}
\end{proof}
\end{proposition}

\begin{remark}
As $\pi_{\GM,n}$ is an \'etale map, differential operators on $\fl$ uniquely lift to $\Sigma_n$. From this one can deduce that $T_\mu^\lambda \omega_{\Sigma_n}^{(a,b)} \isom \pi_{\GM,n}^* T_\mu^\lambda\omega_{\fl}^{(a,b)}$. In other words, translations of the sheaves $\omega_{\Sigma_n}^{(a,b)}$ follow a similar pattern of translations of the sheaves $\omega_{\fl}^{(a,b)}$.
\end{remark}

\subsection{Wall-crossing of equivariant sheaves}
Let $\mu = (a,b)\in\bbZ^2$ be an integral weight. We define the wall-crossing functor as follows:
\begin{align}
    \Theta_s:\Mod_{\tilde{\chi}_\mu}(D_\sigma(G,E))&\to \Mod_{\tilde{\chi}_\mu}(D_\sigma(G,E))\\
    M&\mapsto T_{(-1,0)}^\mu T_\mu^{(-1,0)} M.
\end{align}
In other words, we first translate $M$ to the singular weight, and then we translate back. One also puts $\Theta_s\omega_{\Sigma_n}^{(a,b)} := T_{(-1,0)}^\mu T_\mu^{(-1,0)} \omega_{\Sigma_n}^{(a,b)}$. 
\begin{proposition}\label{prop:1.3.1}
The sheaf $\Theta_s\omega_{\Sigma_n}^{(a,b)}$ is a coherent sheaf on $\Sigma_n$, with $H^i(\Sigma_n,\Theta_s\omega_{\Sigma_n}^{(a,b)})=0$ for $i\ge 1$, and 
\begin{align}
    H^0(\Sigma_n,\Theta_s\omega_{\Sigma_n}^{(a,b)}) \isom \Theta_sH^0(\Sigma_n,\omega_{\Sigma_n}^{(a,b)}).
\end{align}
\begin{proof}
By Proposition \ref{prop:basicoftranslation}, we know that $T_\mu^{(-1,0)}\omega_{\Sigma_n}^{(a,b)}$ has a filtration with graded pieces given by equivariant locally free sheaves of rank $1$, of the form $\omega_{{\Sigma_n}}^{(a',b')}$ for some $(a',b')\in\bbZ^2$. Because the translation functor $T_{(-1,0)}^{\mu}$ is exact, we see that $\Theta_s\omega_{\Sigma_n}^{(a,b)}$ also has a filtration with graded pieces given by equivariant locally free sheaves of rank $1$, of the form $\omega_{{\Sigma_n}}^{(a',b')}$ for some $(a',b')\in\bbZ^2$. From this, one deduces that $\Theta_s\omega_{\Sigma_n}^{(a,b)}$ is a coherent sheaf on $\Sigma_n$, $H^i(\Sigma_n,\Theta_s\omega_{\Sigma_n}^{(a,b)})=0$ for $i\ge 1$, and 
\begin{align}
    H^0(\Sigma_n,\Theta_s\omega_{\Sigma_n}^{(a,b)}) \isom \Theta_sH^0(\Sigma_n,\omega_{\Sigma_n}^{(a,b)}).
\end{align}
\end{proof}
\end{proposition}

First we compute the sheaf $\Theta_s\omega_{\Sigma_n}^{(a,b)}$ with $a\ge b$. Up to twist by determinants, we may assume that $(a,b)=(k,0)$, with $k\ge 0$ an integer. 
\begin{proposition}\label{prop:wcofomegan1}
\begin{enumerate}[(i)]
    \item $T_{(k,0)}^{(-1,0)}\omega_{\Sigma_n}^{(k,0)}\isom \omega_{\Sigma_n}^{(-1,0)}$.
    \item There exists an exact sequence of $\GL_2(L)$-equivariant sheaves on $\Sigma_n$: 
    \begin{align}\label{eq:structureofThetasomegak0}
        0\to \omega_{\Sigma_n}^{(-1,k+1)}\to  \Theta_s\omega_{\Sigma_n}^{(k,0)}\to \omega_{\Sigma_n}^{(k,0)}\to 0.
    \end{align}
\end{enumerate}
\begin{proof}
From the proof of Lemma \ref{lem:basic}, we see that $\omega_{\Sigma_n}^{(k,0)}\ox_E V^{(0,-k-1)}$ is a successive extensions of line bundles $\omega_{\Sigma_n}^{(a-1,-a)}$ with $a=0,1,...,k+1$. Therefore, 
\begin{align}
    T_{(k,0)}^{(-1,0)}\omega_{\Sigma_n}^{(k,0)}= (\omega_{\Sigma_n}^{(k,0)}\ox_E V^{(0,-k-1)})\{\calZ(\frg_\sigma)=\tilde\chi_{(-1,0)}\}\isom \omega_{\Sigma_n}^{(-1,0)}.
\end{align}

The rest is to compute $T_{(-1,0)}^{(k,0)}\omega_{\Sigma_n}^{(-1,0)}$. Using similar methods, we deduce that $\omega_{\Sigma_n}^{(-1,0)}\ox_E V^{(k+1,0)}$ is a successive extension of line bundles $\omega_{\Sigma_n}^{(a-1,k+1-a)}$ with $a=0,1,...,k+1$. After taking $\{\calZ(\frg_\sigma)=\tilde{\chi}_{(k,0)}\}$, we see that only the cases $a=0$ and $a=k+1$ survives. Therefore, the natural maps $\omega_{\Sigma_n}^{(-1,0)}\ox_E V^{(k+1,0)}\surj \omega_{\Sigma_n}^{(-1,0)}\ox_{\calO_{\Sigma_n}}\omega_{\Sigma_n}^{(k+1,0)}=\omega_{\Sigma_n}^{(k,0)}$ and $\omega_{\Sigma_n}^{(-1,k+1)}=\omega_{\Sigma_n}^{(-1,0)}\ox_{\calO_{\Sigma_n}}\omega_{\Sigma_n}^{(0,k+1)}\inj\omega_{\Sigma_n}^{(-1,0)}\ox_E V^{(k+1,0)}$ induce the following exact sequence 
\begin{align}
    0\to \omega_{\Sigma_n}^{(-1,k+1)}\to  \Theta_s\omega_{\Sigma_n}^{(k,0)}\to \omega_{\Sigma_n}^{(k,0)}\to 0.
\end{align}
\end{proof}
\end{proposition}

Next we study the action of $\calZ(\frg_\sigma)$ on $\Theta_s\omega_{\Sigma_n}^{(k,0)}$. By Proposition \ref{prop:wcofomegan1}, we get a $\GL_2(L)$-equivariant exact sequence 
\begin{align}\label{eq:structureofThetasomegak0}
    0\to \omega_{{\Sigma_n}}^{(-1,k+1)}\to  \Theta_s\omega_{{\Sigma_n}}^{(k,0)}\to \omega_{{\Sigma_n}}^{(k,0)}\to 0.
\end{align}
Let $\frc\in\calZ(\frg_\sigma)$ be the Casimir operator and write $c_k = \frac{1}{2}(k+1)^2-\frac{1}{2}$. The action of $\frc-c_k$ on $\omega_{{\Sigma_n}}^{(-1,k+1)}$ and $\omega_{{\Sigma_n}}^{(k,0)}$ is trivial. By considering the action of $\frc$, we get a morphism of short exact sequence:
\begin{center}
\begin{tikzpicture}[descr/.style={fill=white,inner sep=1.5pt}]
    \matrix (m) [
        matrix of math nodes,
        row sep=2.5em,
        column sep=2.5em,
        text height=1.5ex, 
        text depth=0.25ex
    ]
    { 0 & \omega_{{\Sigma_n}}^{(-1,k+1)} & \Theta_s\omega_{{\Sigma_n}}^{(k,0)} & \omega_{{\Sigma_n}}^{(k,0)} & 0  \\
    0 & \omega_{{\Sigma_n}}^{(-1,k+1)} & \Theta_s\omega_{{\Sigma_n}}^{(k,0)} & \omega_{{\Sigma_n}}^{(k,0)} & 0  \\
    };
    \path[->,font=\scriptsize]
    (m-1-1) edge (m-1-2)
    (m-1-2) edge (m-1-3)
    (m-1-3) edge (m-1-4)
    (m-1-4) edge (m-1-5)
    (m-2-1) edge (m-2-2)
    (m-2-2) edge (m-2-3)
    (m-2-3) edge (m-2-4)
    (m-2-4) edge (m-2-5)
    ;
    \path[->,font=\scriptsize]
    (m-1-2) edge node[right] {$\frc-c_k$} (m-2-2)
    (m-1-3) edge node[right] {$\frc-c_k$} (m-2-3)
    (m-1-4) edge node[right] {$\frc-c_k$} (m-2-4)
    ;
\end{tikzpicture}
\end{center}
which is $\GL_2(L)$-equivariant. We denote by $\delta_{\frc}$ the connection map $\omega_{{\Sigma_n}}^{(k,0)}\to \omega_{{\Sigma_n}}^{(-1,k+1)}$ induced by the above morphism of short exact sequence. Besides, starting from the $\GL_2(L)$-equivariant differential map $\calO_{{\Sigma_n}}\to \Omega_{{\Sigma_n}}\isom \omega_{{\Sigma_n}}^{(-1,1)}$, applying $T_{(0,0)}^{(k,0)}$, together with Lemma \ref{lem:T0kOfl}, we get a $\GL_2(L)$-equivariant map $\omega_{{\Sigma_n}}^{(k,0)}\to \omega_{{\Sigma_n}}^{(-1,k+1)}$. We denote this map by $d$. The following proposition is largely inspired by the proof of \cite[6.5.11]{PanII}.
\begin{proposition}\label{prop:ddeltac}
The maps $\delta_{\frc},d:\omega_{{\Sigma_n}}^{(k,0)}\to \omega_{{\Sigma_n}}^{(-1,k+1)}$ coincide up to a non-zero scalar.
\begin{proof}
First we study the case $k=0$. Let $V^{(1,0)}$ be the $\sigma$-standard representation of $\GL_2(L)$ with basis $v_1,v_2\in V$ such that $u^+v_1=0$ and $u^+v_2=v_1$. Let $e_1,e_2\in \omega^{(1,0)}_{{\Sigma_n}}$ be the image of $v_1\ox1, v_2\ox 1$ via $V^{(1,0)}\ox_E\calO_{{\Sigma_n}}\to \omega^{(1,0)}_{{\Sigma_n}}$. Let $U\subset {\Sigma_n}$ be an affinoid open subset. Suppose $e_1$ is invertible on $U$, we have $\omega_{{\Sigma_n}}^{(-1,0)}(U)=\calO_{{\Sigma_n}}(U)e_1^{-1}$. Consider the exact sequence (\ref{eq:structureofThetasomegak0}):
\begin{align}
    0\to \Omega^1_{{\Sigma_n}}\to V^{(1,0)}\ox_E\omega_{{\Sigma_n}}^{(-1,0)}\to \calO_{{\Sigma_n}}\to 0.
\end{align}
Given $s\in \calO_{{\Sigma_n}}(U)$, the element $v_1\ox e_1^{-1}\ox s\in V^{(1,0)}\ox\omega_{{\Sigma_n}}^{(-1,0)}(U)$ maps to $s$. A direct computation shows $h(v_1\ox e_1^{-1}\ox s)=v_1\ox e_1^{-1}\ox hs$, so that $\frac{1}{2}h^2(v_1\ox e_1^{-1}\ox s)=v_1\ox e_1^{-1}\ox \frac{1}{2}h^2s$. As $u^+e_1=0$, we have $u^+(v_1\ox e_1^{-1}\ox s)=v_1\ox e_1^{-1}\ox u^+s$. Then one gets
\begin{align}
    u^-u^+(v_1\ox e_1^{-1}\ox s)=v_2\ox e_1^{-1}\ox u^+s+v_1\ox (-\frac{e_2}{e_1^2})\ox u^+s+v_1\ox e_1^{-1}\ox u^-u^+s.
\end{align}
As $\frc=\frac{1}{2}h^2+h+2u^-u^+$ acts trivially on $s\in \calO_{{\Sigma_n}}(U)$, we see $\frc(v_1\ox e_1^{-1}\ox s)=(v_2\ox e_1^{-1}+e_1\ox -\frac{e_2}{e_1^2})\ox u^+s$. This shows that the connection map $\delta_{\frc}$ satisfies the Lebniz rule: for any $s,t\in \calO_{{\Sigma_n}}(U)$, we have $\delta_{\frc}(st)=\delta_{\frc}(s)t+s\delta_{\frc}(t)$. Moreover, as $v_2\ox e_1^{-1}+v_1\ox -\frac{e_2}{e_1^2}$ maps to zero in via $V^{(1,0)}\ox_E\omega^{(-1,0)}_{{\Sigma_n}}\to \calO_{{\Sigma_n}}$, we see that $v_2\ox e_1^{-1}+v_1\ox -\frac{e_2}{e_1^2}$ is a basis for $\Omega^1_{{\Sigma_n}}$ over $\calO_{{\Sigma_n}}$ on $U$. Besides, on $U$ one directly checks that $d:\calO_{{\Sigma_n}}(U)\to \Omega^1_{{\Sigma_n}}(U)$ is given by $u^+$. From this we see $d$ equals to $\delta_{\frc}$ on $U$ up to a non-zero scalar. By construction both maps are $\GL_2(L)$-equivariant, so that $d$ equals to $\delta_{\frc}$ up to a non-zero scalar.

For general $k$, consider the exact sequence 
\begin{align}
    0\to V^{(k,0)}\ox \Omega^1_{{\Sigma_n}}\to V^{(k,0)}\ox  V^{(1,0)}\ox_E\omega_{{\Sigma_n}}^{(-1,0)}\to V^{(k,0)}\ox  \calO_{{\Sigma_n}}\to 0.
\end{align}
Write $c_k = \frac{1}{2}(k+1)^2-\frac{1}{2}$. The element $\frc-c_k$ induces a morphism of short exact sequence, and induces a connection map $(V^{(k,0)}\ox  \calO_{{\Sigma_n}})[\frc=c_k]\to (V^{(k,0)}\ox \Omega^1_{{\Sigma_n}})[\frc=c_k]$. From the above computation we see this map is exactly the translation $T_{(0,0)}^{(k,0)}$ of $d:\calO_{{\Sigma_n}}\to \Omega^1_{{\Sigma_n}}$ up to a non-zero scalar. This completes the proof.
\end{proof}
\end{proposition}

As a corollary of Proposition \ref{prop:ddeltac}, we deduce the following result.
\begin{corollary}\label{cor:connectionmap}
The action of $\calZ(\frg_\sigma)$ on $\Theta_s\omega_{\Sigma_n}^{(k,0)}$ is not semisimple, and induces a $\GL_2(L)$-equivariant sheaf morphism $\omega_{\Sigma_n}^{(k,0)}\to \omega_{\Sigma_n}^{(-1,k+1)}$. This morphism is given by $d:\omega_{\Sigma_n}^{(k,0)}\to \omega_{\Sigma_n}^{(k,0)}\ox_{\calO_{\Sigma_n}}(\Omega_{\Sigma_n}^{1})^{\ox k+1}$ up to a non-zero scalar.
\end{corollary}

By adjoint properties of translation functors, there exist morphisms of functors $\id\to \Theta_s$ and $\Theta_s\to \id$. In particular, we have natural maps $\omega_{\Sigma_n}^{(k,0)}\to \Theta_s \omega_{\Sigma_n}^{(k,0)}$ and $\Theta_s \omega_{\Sigma_n}^{(k,0)}\to \omega_{\Sigma_n}^{(k,0)}$.


\begin{lemma}\label{lem:WCofd}
The natural map $\Theta_s(\omega_{{\Sigma_n}}^{(k,0)})\to \omega_{{\Sigma_n}}^{(k,0)}$ is surjective. The natural map $\omega_{{\Sigma_n}}^{(-1,k+1)}\to \Theta_s(\omega_{{\Sigma_n}}^{(-1,k+1)})$ is injective.
\begin{proof}
We first handle the case $k=0$ and the general case follows by applying $T_{(0,0)}^{(k,0)}$. Let $U\subset{\Sigma_n}$ be an open affinoid subset such that $e_1$ vanishes nowhere on $U$. We have the following exact sequence 
\begin{align}
    0\to \omega^{(-1,0)}_{{\Sigma_n}}\to (V^{(1,0)})^*\ox\calO_{{\Sigma_n}}\to \omega_{{\Sigma_n}}^{(0,-1)}\to 0.
\end{align}
Let $v_1^*,v_2^*$ be basis vectors of $(V^{(1,0)})^*$, and let $e_1^*,e_2^*$ be the image of $v_1^*\ox 1,v_2^*\ox 1$ under the map $(V^{(1,0)})^*\ox\calO_{{\Sigma_n}}\to \omega_{{\Sigma_n}}^{(0,-1)}$. In this case, $e_2^*$ is invertible on $U$, and $v_1^*\ox 1+v_2^*\ox-\frac{e_1^*}{e_2^*}$ is a basis for $\omega^{(-1,0)}_{{\Sigma_n}}\subset (V^{(1,0)})^*\ox\calO_{{\Sigma_n}}$. Then we see $v_1\ox (v_1^*\ox 1+v_2^*\ox-\frac{e_1^*}{e_2^*})\ox s\in V^{(1,0)}\ox \omega^{(-1,0)}_{{\Sigma_n}}$ is mapped to any given $s\in\calO_{{\Sigma_n}}(U)$. As $\Theta_s(\calO_{{\Sigma_n}})\to \calO_{{\Sigma_n}}$ is $\GL_2(L)$-equivariant, this shows $\Theta_s(\calO_{{\Sigma_n}})(U)\to \calO_{{\Sigma_n}}(U)$ is surjective for any open affinoid subset $U\subset{\Sigma_n}$.

Next, we show that $\Omega^1_{{\Sigma_n}}\to \Theta_s(\Omega^1_{{\Sigma_n}})$ is injective. Let $U\subset{\Sigma_n}$ be an open affinoid subset such that $e_1$ does not vanish. We have the following exact sequence
\begin{align}
    0\to \omega^{(-2,1)}_{{\Sigma_n}}\to (V^{(1,0)})^*\ox\Omega^1_{{\Sigma_n}}\to \omega_{{\Sigma_n}}^{(-1,0)}\to 0.
\end{align}
Let $dx$ be a basis of $\Omega^1_{{\Sigma_n}}(U)$ over $\calO_{{\Sigma_n}}(U)$. Then $v_1^*\ox dx,v_2^*\ox dx$ are respectively mapped to $e_1^*dx,e_2^*dx$. Now given any $sdx\in \Omega^1_{{\Sigma_n}}(U)=\calO_{{\Sigma_n}}(U)dx$ with $s\in \calO_{{\Sigma_n}}(U)$, we see $s\mapsto v_1\ox v_1^*\ox sdx+v_2\ox v_2^*\ox sdx\in V^{(1,0)}\ox (V^{(1,0)})^*\ox \Omega^1_{{\Sigma_n}}$ which identifies $\Omega^1_{{\Sigma_n}}$ with a subsheaf of $\Theta_s(\Omega^1_{{\Sigma_n}})$. This shows that $\Omega^1_{{\Sigma_n}}\to \Theta_s(\Omega^1_{{\Sigma_n}})$ is injective. 
\end{proof}
\end{lemma}

The following lemma relates the wall-crossing of $\omega_{\Sigma_n}^{(k,0)}$ with the wall-crossing of $\omega_{\Sigma_n}^{(-1,k+1)}$.
\begin{lemma}\label{lem:Thetasdisom}
For any $k\ge 0$, the wall-crossing of the differential operator $\Theta_s(d):\Theta_s(\omega_{{\Sigma_n}}^{(k,0)})\to \Theta_s(\omega_{{\Sigma_n}}^{(-1,k+1)})$ is an isomorphism. 
\begin{proof}
We first handle the case $k=0$ and the general case follows from applying $T_{(0,0)}^{(k,0)}$. We will show that $T_{(0,0)}^{(-1,0)}d:T_{(0,0)}^{(-1,0)}\calO_{{\Sigma_n}}\to T_{(0,0)}^{(-1,0)}\Omega^1_{{\Sigma_n}}$ is an isomorphism. Consider the following diagram with exact rows:
\begin{center}
\begin{tikzpicture}[descr/.style={fill=white,inner sep=1.5pt}]
    \matrix (m) [
        matrix of math nodes,
        row sep=2.5em,
        column sep=2.5em,
        text height=1.5ex, 
        text depth=0.25ex
    ]
    { 0 & \omega^{(-1,0)}_{{\Sigma_n}} & (V^{(1,0)})^*\ox \calO_{{\Sigma_n}} & \omega^{(0,-1)}_{{\Sigma_n}} & 0  \\
    0 & \omega^{(-2,1)}_{{\Sigma_n}} & (V^{(1,0)})^*\ox \Omega^1_{{\Sigma_n}} & \omega^{(-1,0)}_{{\Sigma_n}} & 0  \\
    };
    \path[->,font=\scriptsize]
    (m-1-1) edge (m-1-2)
    (m-1-2) edge (m-1-3)
    (m-1-3) edge (m-1-4)
    (m-1-4) edge (m-1-5)
    (m-2-1) edge (m-2-2)
    (m-2-2) edge (m-2-3)
    (m-2-3) edge (m-2-4)
    (m-2-4) edge (m-2-5)
    ;
    \path[->,font=\scriptsize]
    (m-1-3) edge node[right] {$1\ox d$} (m-2-3)
    ;
\end{tikzpicture}
\end{center}
Let $U\subset{\Sigma_n}$ be an affinoid open subset such that $e_1$ does not vanish on $U$. Given $s\in \calO_{{\Sigma_n}}(U)$, then $v_1^*\ox 1\ox s+v_2^*\ox-\frac{e_1^*}{e_2^*}\ox s\in \omega^{(-1,0)}_{{\Sigma_n}} \subset (V^{(1,0)})^*\ox \calO_{{\Sigma_n}}$ is mapped to $v_1^*\ox ds+v_2^*\ox (-d\frac{e_1^*}{e_2^*}s)+v_2^*\ox -\frac{e_1^*}{e_2^*}ds$ under $1\ox d$, so that it is mapped to $-e_2^*d\frac{e_1^*}{e_2^*}s$. Note that under the identification $V^{(1,0)}\isom (V^{(1,0)})^*\ox \det$, $e_1$ is identified with $e_2^*\ox \det$ so that $\frac{e_1^*}{e_2^*}$ is identified with $\frac{e_2}{e_1}$. In particular, $d\frac{e_1^*}{e_2^*}$ is a basis of $\Omega^1_{{\Sigma_n}}(U)$ over $\calO_{{\Sigma_n}}(U)$, so $e_2^*d\frac{e_1^*}{e_2^*}$ is a basis of $\omega^{(-1,0)}_{{\Sigma_n}}=\omega^{(0,-1)}_{{\Sigma_n}}\ox_{\calO_{{\Sigma_n}}} \Omega^1_{{\Sigma_n}}$. This shows that $d$ induces an isomorphism $T_{(0,0)}^{(-1,0)}d:T_{(0,0)}^{(-1,0)}\calO_{{\Sigma_n}}\to T_{(0,0)}^{(-1,0)}\Omega^1_{{\Sigma_n}}$. Therefore, $d$ induces an isomorphism $\Theta_s(d):\Theta_s\calO_{{\Sigma_n}}\to \Theta_s\Omega^1_{{\Sigma_n}}$.
\end{proof}
\end{lemma}

\begin{corollary}
There exists an exact sequence of $\GL_2(L)$-equivariant sheaves on $\Sigma_n$: 
\begin{align}\label{eq:structureofThetasomegak01}
    0\to \omega_{\Sigma_n}^{(-1,k+1)}\to  \Theta_s\omega_{\Sigma_n}^{(-1,k+1)}\to \omega_{\Sigma_n}^{(k,0)}\to 0.
\end{align}
\begin{proof}
By Lemma \ref{lem:Thetasdisom}, we know that $\Theta_s(d):\Theta_s(\omega_{\Sigma_n}^{(k,0)})\to \Theta_s(\omega_{\Sigma_n}^{(-1,k+1)})$ is an isomorphism. From this, together with (\ref{eq:structureofThetasomegak0}), we deduce the exact sequence (\ref{eq:structureofThetasomegak01}).
\end{proof}
\end{corollary}

From Lemma \ref{lem:Thetasdisom}, together with Lemma \ref{lem:WCofd} we deduce the following corollary.
\begin{corollary}\label{cor:ffffff}
For any $k\ge 0$, the natural map $\omega^{(k,0)}_{{\Sigma_n}}\to \Theta_s\omega^{(k,0)}_{{\Sigma_n}}$ is given by the composition $d:\omega^{(k,0)}_{{\Sigma_n}}\to \omega^{(-1,k+1)}_{{\Sigma_n}}$ with the natural inclusion $\omega^{(-1,k+1)}_{{\Sigma_n}}\subset \Theta_s\omega^{(-1,k+1)}_{{\Sigma_n}}\isom \Theta_s\omega^{(k,0)}_{{\Sigma_n}}$. The natural map $\Theta_s\omega^{(-1,k+1)}_{{\Sigma_n}}\to \omega^{(-1,k+1)}_{{\Sigma_n}}$ is given by the composition $\Theta_s\omega^{(-1,k+1)}_{{\Sigma_n}}\isom \Theta_s\omega^{(k,0)}_{{\Sigma_n}}\to \omega^{(k,0)}_{{\Sigma_n}}$ with $d:\omega^{(k,0)}_{{\Sigma_n}}\to \omega^{(-1,k+1)}_{{\Sigma_n}}$.
\end{corollary}

\subsection{Wall-crossing of some locally analytic $\GL_2(L)$-representations}
In this section, we apply the results on wall-crossing of $\omega_{\Sigma_n}^{(a,b)}$ to compute the wall-crossing of some locally analytic representations of $\GL_2(L)$ defined using the cohomology of these sheaves. 

Let $\tau$ be a smooth irreducible representation of $D_L^\times$ over $E$, such that the action of $D_L^\times$ does not factor through the reduced norm map, and $\varpi\in L^\times\subset D_L^\times$ acts trivially on $\tau$. Fix an integer $n\ge 1$, and we assume that $1+\varpi^n\calO_{D_L}$ acts trivially on $\tau$. As $\Sigma_n$ carries an action of $D_L^\times$ such that $1+\varpi^n\calO_{D_L}$ acts trivially, we can cut out some locally analytic representations of $\GL_2(L)$ using $\tau$.

Starting from the de Rham complex $\calO_{\Sigma_n}\to \Omega^1_{\Sigma_n}$ of $\Sigma_n$, taking global sections, we get an exact sequence 
\begin{align}
    0\to H^0_{\dR}(\Sigma_n)\to H^0(\Sigma_n,\calO_{\Sigma_n})\to H^0(\Sigma_n,\Omega_{\Sigma_n}^1)\to H^1_{\dR}(\Sigma_n)\to 0. \label{eq:dR0}
\end{align}
where $H^i_{\dR}(\Sigma_n)$ is the $i$th de Rham cohomology of $\Sigma_n$, for $i=0,1$. In \cite[\S 4]{DLB17}, \cite[Th\'eor\'eme 5.8]{CDN20}, they give a detailed description of $H^1_{\dR}(\Sigma_{n,C})$ where $\Sigma_{n,C}:=\Sigma_n\times_E C$. For $i=0,1$, we have a natural isomorphism 
\begin{align}
    H^i_{\dR}(\Sigma_n)\hat\ox_E C\isom H^i_{\dR}(\Sigma_{n,C}).
\end{align}
Indeed, for $\{U_i\}$ a Stein cover of $\Sigma_n$, $\{U_i\times_E C\}$ is a Stein cover of $\Sigma_{n,C}$. Therefore, 
\begin{align}
    \calO(\calM_n^{\varpi})=\ilim_i\calO(U_i\times_E C)\isom \ilim_i (\calO(U_i)\hat\ox_E C)\ov{(*)}\isom (\ilim_i\calO(U_i))\hat\ox_E C\isom \calO(\Sigma_n)\hat\ox_E C
\end{align}
where $(*)$ follows by \cite[Proposition 1.1.29]{Eme17}. By a similar proof, $\Omega^1(\Sigma_{n,C})\isom \Omega^1(\Sigma_n)\hat\ox_E C$. From this we deduce that $H^1_{\dR}(\Sigma_n)\hat\ox_E C\isom H^1_{\dR}(\Sigma_{n,C})$.

Let $\pi_{\infty ,C}$ be a smooth supercuspidal representation of $\GL_2(L)$ over $C$, corresponds to $\tau\hat\ox_E C$ via the Jacquet-Langlands correspondence. From the classification of smooth supercuspidal representations of $\GL_2(L)$, we know that $\pi_{\infty ,C}$ is canonically defined over $\bar L$. Fix a choice of $q^{1/2}$ in $\bar L$. By \cite[Chapter 4]{BS07}, after enlarge $E$ if necessary, there exists a unique smooth supercuspidal representation of $\GL_2(L)$ over $E$, which we denote it by $\pi_{\infty }$, such that $\pi_{\infty }\hat\ox_E C\isom \pi_{\infty ,C}$. We equip $\pi_\infty$ with the ind-finite-dimensional topology. Then 
\begin{align}
    \Hom_{D_L^\times}(\tau,H^1_{\dR}(\Sigma_n))\isom (\pi_{\infty }')^{\oplus 2}.
\end{align}
where $\pi_{\infty }'$ is the strong dual of $\pi_{\infty }$. Indeed, by \cite[Th\'eor\'eme 5.8]{CDN20}, as $1+\varpi^n\calO_{D_L}$ acts trivially on $\tau$, we know that $\GL_2(L)$-equivariantly
\begin{align}
    \Hom_{D_L^\times}(\tau\hat\ox_E C,H^1_{\dR}(\calM_n^{\varpi}))\isom (\pi_{\infty ,C}')^{\oplus 2}
\end{align}
where $\pi_{\infty ,C}$ is a smooth supercuspidal representation of $\GL_2(L)$ over $C$, and $\pi_{\infty ,C}'$ is the strong dual of $\pi_{\infty ,C}$.  Therefore, 
\begin{align}
    \Hom_{\GL_2(L)}(\pi_{\infty }',\Hom_{D_L^\times}(\tau,H^1_{\dR}(\Sigma_n)))\hat\ox_E C\isom \Hom_{\GL_2(L)}(\pi_{\infty ,C}',(\pi_{\infty ,C}')^{\oplus 2})
\end{align}
is $2$-dimensional. Hence $\Hom_{\GL_2(L)}(\pi_{\infty }',\Hom_{D_L^\times}(\tau,H^1_{\dR}(\Sigma_n)))$ is $2$-dimensional. Similarly one shows that $\Hom_{D_L^\times}(\tau,H^1_{\dR}(\Sigma_n))$ has no other components. From this we deduce that
\begin{align}
    \Hom_{D_L^\times}(\tau,H^1_{\dR}(\Sigma_n))\isom (\pi_{\infty }')^{\oplus 2}.
\end{align}

Applying the translation functor $T_{(0,0)}^{(k,0)}$ for $k\ge 1$ to (\ref{eq:dR0}), we get an exact sequence 
\begin{align}
    0\to V^{(k,0)}\ox_E H^0_{\dR}(\Sigma_n)\to H^0(\Sigma_n,\omega^{(k,0)}_{\Sigma_n})\to H^0(\Sigma_n,\omega_{\Sigma_n}^{(-1,k+1)})\to V^{(k,0)}\ox_E H^1_{\dR}(\Sigma_n)\to 0 \label{eq:dRk}
\end{align}
Indeed, by Proposition \ref{prop:basicoftranslation}, we know that translation functor commutes with taking cohomology. Besides, by \cite[Th\'eor\`eme 4.1]{CDN20}, we know that $\frg_\sigma$ acts trivially on $H^i_{\dR}(\Sigma_n)$ for $i=0,1$. Hence by definition of the translation functor we deduce $T_{(0,0)}^{(k,0)}H^i_{\dR}(\Sigma_n)\isom V^{(k,0)}\ox_E H^i_{\dR}(\Sigma_n)$. Note that $\varpi\in D_L^\times$ acts trivially on (\ref{eq:dRk}). Therefore, applying $\Hom_{D_L^\times}(\tau,-)=\Hom_{{\varpi^{\bbZ}}}(\tau,-)^{{D_L^\times/{\varpi^{\bbZ}}}}$ to (\ref{eq:dRk}), as $D_L^\times/{\varpi^{\bbZ}}$ is compact, we get an exact sequence
\begin{align}
    0\to V^{(k,0)}\ox_E \Hom_{D_L^\times}(\tau,H^0_{\dR}(\Sigma_n))&\to \Hom_{D_L^\times}(\tau,H^0(\Sigma_n,\omega^{(k,0)}_{\Sigma_n}))\\&\to \Hom_{D_L^\times}(\tau,H^0(\Sigma_n,\omega_{\Sigma_n}^{(-1,k+1)}))\to V^{(k,0)}\ox\Hom_{D_L^\times}(\tau,H^1_{\dR}(\Sigma_n))\to 0.
\end{align} 
Besides, the action of $D_L^\times$ on $H^0_{\dR}(\Sigma_n)$ factors through the reduced norm. This implies $\Hom_{D_L^\times}(\tau,H^0_{\dR}(\Sigma_n))=0$. Consequently, we get an exact sequence 
\begin{align}
    0\to \Hom_{D_L^\times}(\tau,H^0(\Sigma_n,\omega^{(k,0)}_{\Sigma_n}))\to \Hom_{D_L^\times}(\tau,H^0(\Sigma_n,\omega_{\Sigma_n}^{(-1,k+1)}))\to V^{(k,0)}\ox_E\Hom_{D_L^\times}(\tau,H^1_{\dR}(\Sigma_n))\to 0.
\end{align}

To simplify the notation, we name each term in the exact sequence as follows: 
\begin{itemize}
    \item $M_{\alg,k}:=V^{(k,0)}\ox_E\pi_{\infty }'$.
    \item $\tilde{M}_k:=\Hom_{E[D_L^\times]}(\tau,H^0(\Sigma_n,\omega_{\Sigma_n}^{(-1,k+1)}))$.
    \item $M_{c,k}:=\Hom_{E[D_L^\times]}(\tau,H^0(\Sigma_n,\omega^{(k,0)}_{\Sigma_n}))$.
\end{itemize}

Then we have an exact sequence of $D_\sigma(G,E)$-modules 
\begin{align}
    0\to M_{c,k}\to \tilde{M}_k\to M_{\alg,k}^{\oplus 2}\to 0. \label{eq:Mtilde}
\end{align}

We first compute the wall-crossing of each $D_\sigma(G,E)$-modules above.
\begin{proposition}\label{prop:1.4.3}
$\Theta_s M_{\alg,k} = 0$, and the natural map $M_{c,k}\to \tilde{M}_k$ induces an isomorphism $\Theta_s(M_{c,k})\aisom \Theta_s(\tilde{M}_k)$.
\begin{proof}
We first compute $T_{(k,0)}^{(-1,0)}M_{\alg,k}$. By definition, we have 
\begin{align}
    T_{(k,0)}^{(-1,0)}M_{\alg,k} = (M_{\alg,k}\ox_E V^{(0,-k-1)})\{\calZ(\frg_\sigma)=\tilde{\chi}_{(-1,0)}\}&=(\pi_\infty'\ox_E V^{(k,0)}\ox_E V^{(0,-k-1)})\{\calZ(\frg_\sigma)=\tilde{\chi}_{(-1,0)}\}
\end{align} 
As the action of $\frg_\sigma$ on $\pi_\infty'$ is trivial, we have 
\begin{align}
    (\pi_\infty'\ox_E V^{(k,0)}\ox_E V^{(0,-k-1)})\{\calZ(\frg_\sigma)=\tilde{\chi}_{(-1,0)}\}=\pi_\infty'\ox_E (V^{(k,0)}\ox_E V^{(0,-k-1)})\{\calZ(\frg_\sigma)=\tilde{\chi}_{(-1,0)}\}.
\end{align}
However, $V^{(k,0)}\ox_E V^{(0,-k-1)}$ is an algebraic representation of $\GL_2(L)$. Therefore, $(V^{(k,0)}\ox_E V^{(0,-k-1)})\{\calZ(\frg_\sigma)=\tilde{\chi}_{(-1,0)}\}=0$. Finally, as $\Theta_s$ is an exact functor, applying $\Theta_s(-)$ to (\ref{eq:Mtilde}), we have an exact sequence 
\begin{align}
    0\to \Theta_s(M_{c,k})\to \Theta_s(\tilde{M}_k)\to \Theta_s(M_{\alg,k}^{\oplus 2})\to 0.
\end{align}
Since $\Theta_s(M_{\alg,k})=0$, we deduce that $\Theta_s(M_{c,k})\isom \Theta_s(\tilde{M}_k)$.
\end{proof}
\end{proposition}

\begin{proposition}\label{prop:injsurj}
The wall-crossing of $M_{c,k}$ fits into an exact sequence 
\begin{align}
    0\to \tilde{M}_k\to \Theta_s(M_{c,k})\to M_{c,k}\to 0.
\end{align}
The natural map $M_{c,k}\to \Theta_s(M_{c,k})$ is injective, and the natural map $\Theta_s(M_{c,k})\to M_{c,k}$ is surjective. 
\begin{proof}
By Proposition \ref{prop:wcofomegan1}, we get an exact sequence 
\begin{align}
    0\to \omega_{\Sigma_n}^{(-1,k+1)}\to  \Theta_s\omega_{\Sigma_n}^{(k,0)}\to \omega_{\Sigma_n}^{(k,0)}\to 0.
\end{align}
Taking $H^0(\Sigma_n,-)$, we get an exact sequence 
\begin{align}
    0\to H^0(\Sigma_n,\omega_{\Sigma_n}^{(-1,k+1)})\to  \Theta_sH^0(\Sigma_n,\omega_{\Sigma_n}^{(k,0)})\to H^0(\Sigma_n,\omega_{\Sigma_n}^{(k,0)})\to 0
\end{align}
because $H^1(\Sigma_n,\omega_{\Sigma_n}^{(-1,k+1)})=0$, and $\Theta_s$ commutes with taking cohomology, see Proposition \ref{prop:1.3.1}. Applying $\Hom_{D_L^\times}(\tau,-)$ (this functor commutes with $\Theta_s(-)$), we get an exact sequence 
\begin{align}
    0\to \tilde{M}_k\to \Theta_s(M_{c,k})\to M_{c,k}\to 0
\end{align}
because $\varpi$ acts trivially on $\tau$ and $\Sigma_n$ and $D_L^\times/\varpi^{\bbZ}$ is compact. Similarly, one can show that $M_{c,k}\to \Theta_s(M_{c,k})$ is injective.
\end{proof}
\end{proposition}

\begin{corollary}\label{cor:nonzero}
The $D_\sigma(G,E)$-module $M_{c,k}$ is non-zero.
\begin{proof}
Suppose $M_{c,k}= 0$. Then by Proposition \ref{prop:injsurj}, we see $\tilde{M}_k=0$. But $\tilde{M}_k\neq 0$, because it has a non-zero quotient $M_{\alg,k}^{\oplus 2}$, see (\ref{eq:Mtilde}).
\end{proof}
\end{corollary}

\subsection{Calculations of $\Ext^1_{D_\sigma(G,E)}(M_{\alg,k},M_{c,k})$}
In this section, we use the technique of wall-crossing to compute the extension groups $\Ext^1_{D_\sigma(G,E)}(M_{c,k},M_{\alg,k})$. We first study morphisms between these representations. 
\begin{lemma}\label{lem:endoofalg}
We have $\Hom_{D_\sigma(G,E)}(M_{\alg,k},M_{\alg,k})=E$.
\begin{proof}
Using adjoint properties of translation functors, we may reduce to the case $k=0$. Then $M_{\alg,k}$ is the strong dual of smooth irreducible admissible representation of $\GL_2(L)$. Therefore, $\Hom_{D_\sigma(G,E)}(M_{\alg,k},M_{\alg,k})=E$ by \cite[Corollary 3.3]{Sch02Dual}.
\end{proof}
\end{lemma}

\begin{proposition}\label{prop:Hom}
We have $\Hom_{D_\sigma(G,E)}(M_{c,k},M_{\alg,k})=0$ and $\Hom_{D_\sigma(G,E)}(M_{\alg,k},M_{c,k})=0$.
\begin{proof}
Let $\phi:M_{c,k}\to M_{\alg,k}$ be any $\GL_2(L)$-equivariant morphism. By Proposition \ref{prop:1.4.3}, the wall-crossing of $M_{\alg,k}$ is $0$, so that the wall-crossing of $\phi$ gives a commutative diagram 
\begin{center}
    \begin{tikzpicture}[descr/.style={fill=white,inner sep=1.5pt}]
        \matrix (m) [
            matrix of math nodes,
            row sep=2.5em,
            column sep=2.5em,
            text height=1.5ex, 
            text depth=0.25ex
        ]
        { \Theta_sM_{c,k}&0  \\
          M_{c,k}&M_{\alg,k} \\
        };
        \path[->,font=\scriptsize]
        (m-1-1) edge node[auto] {$\Theta_s\phi$} (m-1-2)
        (m-1-2) edge node[auto] {} (m-2-2)
        (m-2-1) edge node[auto] {$\phi$} (m-2-2)
        ;
        \path[->>,font=\scriptsize]
        (m-1-1) edge node[auto] {} (m-2-1)
        ;
    \end{tikzpicture}.
\end{center}
By Proposition \ref{prop:injsurj}, the natural map $\Theta_sM_{c,k}\to M_{c,k}$ is surjective. From this, one deduces that $\phi=0$. Next, let $\phi:M_{\alg,k}\to M_{c,k}$ be any $\GL_2(L)$-equivalent morphism. The wall-crossing of $\phi$ gives a commutative diagram 
\begin{center}
    \begin{tikzpicture}[descr/.style={fill=white,inner sep=1.5pt}]
        \matrix (m) [
            matrix of math nodes,
            row sep=2.5em,
            column sep=2.5em,
            text height=1.5ex, 
            text depth=0.25ex
        ]
        { M_{\alg,k}&M_{c,k}  \\
          0&\Theta_sM_{c,k} \\
        };
        \path[->,font=\scriptsize]
        (m-1-1) edge node[auto] {$\phi$} (m-1-2)
        (m-1-1) edge node[auto] {} (m-2-1)
        (m-2-1) edge node[auto] {$\Theta_s\phi$} (m-2-2)
        ;
        \path[right hook->]
        (m-1-2) edge node[auto] {} (m-2-2)
        ;
    \end{tikzpicture}.
\end{center}
By Proposition \ref{prop:injsurj}, the natural map $ M_{c,k} \to \Theta_s M_{c,k}$ is injective. From this one deduces that $\phi=0$.
\end{proof}
\end{proposition}

\begin{corollary}\label{cor:semisimpleonext}
Let $M$ be a $D_\sigma(G,E)$-module, which is an extension of $M_{c,k}$ by $M_{\alg,k}$. Then the action of $\calZ(\frg_\sigma)$ on $M$ is semisimple. 

Similarly, let $W$ be a $D_\sigma(G,E)$-module, which is an extension of $M_{\alg,k}$ by $M_{c,k}$. Then the action of $\calZ(\frg_\sigma)$ on $W$ is semisimple.
\begin{proof}
Let $M$ be a $D_\sigma(G,E)$-module, which is an extension of $M_{c,k}$ by $M_{\alg,k}$. Let $\frc\in\calZ(\frg_\sigma)$ be the Casimir operator. Consider the following morphism of exact sequences of $D_\sigma(G,E)$-modules:
\begin{center}
\begin{tikzpicture}[descr/.style={fill=white,inner sep=1.5pt}]
    \matrix (m) [
        matrix of math nodes,
        row sep=2.5em,
        column sep=2.5em,
        text height=1.5ex, 
        text depth=0.25ex
    ]
    { 0 & M_{\alg,k} & M & M_{c,k} & 0  \\
    0 & M_{\alg,k} & M & M_{c,k} & 0  \\
    };
    \path[->,font=\scriptsize]
    (m-1-1) edge (m-1-2)
    (m-1-2) edge (m-1-3)
    (m-1-3) edge (m-1-4)
    (m-1-4) edge (m-1-5)
    (m-2-1) edge (m-2-2)
    (m-2-2) edge (m-2-3)
    (m-2-3) edge (m-2-4)
    (m-2-4) edge (m-2-5)
    ;
    \path[->,font=\scriptsize]
    (m-1-2) edge node[right] {$\frc-c_k$} (m-2-2)
    (m-1-3) edge node[right] {$\frc-c_k$} (m-2-3)
    (m-1-4) edge node[right] {$\frc-c_k$} (m-2-4)
    ;
\end{tikzpicture}
\end{center}
with $c_k = \frac{1}{2}(k+1)^2-\frac{1}{2}$. This induces a connection map of $D_\sigma(G,E)$-modules $M_{c,k}\to M_{\alg,k}$. By Proposition \ref{prop:Hom}, we know this connection map is zero. Hence $M=M[\frc=c_k]$. Similarly we have $M=M[z=k]$ where $z=\left( \begin{matrix}1&0\\0&1\end{matrix} \right)\in\calZ(\frg_\sigma)$. Therefore, $M$ has infinitesimal character $\tilde{\chi}_{(k,0)}$. Using similar methods we deduce results for extensions of $M_{\alg,k}$ by $M_{c,k}$.
\end{proof}
\end{corollary}

Recall that by Proposition \ref{prop:injsurj}, the natural maps $M_{c,k}\to \Theta_s(M_{c,k})$ is injective and $\Theta_s(M_{c,k})\to M_{c,k}$ is surjective. Our next result show that if we replace $M_{c,k}$ by extensions of $M_{c,k}$ and $M_{\alg,c}$, similar results still holds.
\begin{theorem}\label{thm:quo}
Let $W_k$ be a $D_\sigma(G,E)$-module, which is a non-split extension of $M_{c,k}$ by $M_{\alg,k}$. Then $\Theta_s(W_k)\isom \Theta_s(M_{c,k})$, and the natural map $\Theta_s(W_k)\to W_k$ is surjective.
\begin{proof}
First of all, as $\Theta_s$ is exact and $\Theta_s(M_{\alg,k})=0$ by Proposition \ref{prop:1.4.3}, we deduce that the natural map $W_k\to M_{c,k}$ induces an isomorphism $\Theta_s(W_k)\isom \Theta_s(M_{c,k})$. The wall-crossing of the map $W_k\surj M_{c,k}$ induces a commutative diagram 
\begin{center}
\begin{tikzpicture}[descr/.style={fill=white,inner sep=1.5pt}]
    \matrix (m) [
        matrix of math nodes,
        row sep=2.5em,
        column sep=2.5em,
        text height=1.5ex, 
        text depth=0.25ex
    ]
    { \Theta_s W_k&\Theta_s M_{c,k}  \\
      W_k& M_{c,k} \\
    };
    \path[->,font=\scriptsize]
    (m-1-1) edge node[auto] {$\sim$} (m-1-2)
    (m-1-1) edge node[auto] {$\alpha$} (m-2-1)
    ;
    \path[->>,font=\scriptsize]
    (m-2-1) edge node[auto] {$\beta$} (m-2-2)
    (m-1-2) edge node[auto] {} (m-2-2)
    ;
\end{tikzpicture}.

\end{center}
By the kernel-cokernel lemma, we observe that $\coker \alpha$ is a quotient of $\ker \beta = M_{\alg,k}$ as $\coker \beta\comp \alpha = 0$. Since $M_{\alg,k}$ is irreducible, $\coker \alpha$ is either isomorphic to $M_{\alg,k}$ or $0$. If $\coker \alpha\isom M_{\alg,k}$, then $M_{\alg,k}$ is a quotient of $W_{k}$. However, 
\begin{align}
    \Hom_{D_\sigma(G,E)}(W_k,M_{\alg,k})\inj\Hom_{D_\sigma(G,E)}(M_{\alg,k},M_{\alg,k})
\end{align}
as $\Hom_{D_\sigma(G,E)}(M_{c,k},M_{\alg,k})=0$ by Proposition \ref{prop:Hom}. This shows that $W_k\to \coker\alpha\isom M_{\alg,k}$ comes from a lift of an endomorphism of $M_{\alg,k}$. By Lemma \ref{lem:endoofalg}, we see this lifting gives a splitting of the extension $W_k\isom M_{\alg,k}-M_{c,k}$. This is a contradiction. Therefore, $\coker \alpha=0$, which implies that $\alpha$ is surjective.
\end{proof}
\end{theorem}

\begin{lemma}\label{lem:wcsemisimple}
We have $\Theta_s(M_{c,k})[\calZ(\frg_\sigma)=\tilde{\chi}_{(k,0)}]\isom \tilde{M}_k$. Similarly, let $\tilde W_k$ be the maximal quotient of $\Theta_s(M_{c,k})$ such that $\calZ(\frg_\sigma)$ acts via $\tilde{\chi}_{(k,0)}$. Then $\tilde W_k$ is an extension of $M_{c,k}$ by $M_{\alg,k}^{\oplus 2}$.
\begin{proof}
Let $\frc\in\calZ(\frg_\sigma)$ be the Casimir operator, and $c_k=\frac{1}{2}(k+1)^2-\frac{1}{2}$. By Proposition \ref{prop:injsurj}, we may consider the following morphism of exact sequences of $D_\sigma(G,E)$-modules 
\begin{center}
\begin{tikzpicture}[descr/.style={fill=white,inner sep=1.5pt}]
    \matrix (m) [
        matrix of math nodes,
        row sep=2.5em,
        column sep=2.5em,
        text height=1.5ex, 
        text depth=0.25ex
    ]
    { 0 & \tilde{M}_k & \Theta_s (M_{c,k}) & M_{c,k} & 0  \\
    0 & \tilde{M}_k & \Theta_s(M_{c,k}) & M_{c,k} & 0  \\
    };
    \path[->,font=\scriptsize]
    (m-1-1) edge (m-1-2)
    (m-1-2) edge (m-1-3)
    (m-1-3) edge (m-1-4)
    (m-1-4) edge (m-1-5)
    (m-2-1) edge (m-2-2)
    (m-2-2) edge (m-2-3)
    (m-2-3) edge (m-2-4)
    (m-2-4) edge (m-2-5)
    ;
    \path[->,font=\scriptsize]
    (m-1-2) edge node[right] {$\frc-c_k$} (m-2-2)
    (m-1-3) edge node[right] {$\frc-c_k$} (m-2-3)
    (m-1-4) edge node[right] {$\frc-c_k$} (m-2-4)
    ;
\end{tikzpicture}.
\end{center}
By Corollary \ref{cor:connectionmap}, we know the connection map is induced by $d:\omega_{\Sigma_n}^{(k,0)}\to \omega_{\Sigma_n}^{(-1,k+1)}$, which induces an injection $M_{c,k}\inj \tilde{M}_k$ (here we use that $\tau$ does not appear in $H^0_{\dR}(\Sigma_n)$). From this we deduce $\Theta_s(M_{c,k})[\calZ(\frg_\sigma)=\tilde{\chi}_{(k,0)}]\isom \tilde{M}_k$ and $\tilde W_k$ fits into an exact sequence 
\begin{align}
    0\to M_{\alg,k}^{\oplus 2}\to \tilde{W}_k\to M_{c,k}\to 0
\end{align}
\end{proof}
\end{lemma}

\begin{corollary}\label{cor:quo}
Let $W_k$ be a $D_\sigma(G,E)$-module, which is an extension of $M_{c,k}$ by $M_{\alg,k}$. The natural map $\Theta_s(W_k)\to W_k$ factors through the quotient $\Theta_s(W_k)\surj \tilde{W}_k$ and induces a surjection $\tilde{W}_k\surj W_k$.
\begin{proof}
As the action of $\calZ(\frg_\sigma)$ on $W_k$ is semisimple by Corollary \ref{cor:semisimpleonext}, the map $\Theta_s(M_{c,k})\isom \Theta_s(W_k)\to W_k$ factors through the maximal semisimple quotient of $\Theta_s(W_k)$ such that $\calZ(\frg_\sigma)$ acts via $\tilde{\chi}_{(k,0)}$, by Lemma \ref{lem:wcsemisimple}. This induces a map $\tilde{W}_k\surj W_k$.
\end{proof}
\end{corollary}

Similarly, we have 
\begin{theorem}\label{thm:inj}
Let $M_k$ be a $D_\sigma(G,E)$-module, which is a non-split extension of $M_{\alg,k}$ by $M_{c,k}$. Then $\Theta_s(M_k)\isom \Theta_s(M_{c,k})$, and the natural map $M_k\to \Theta_s(M_k)$ is injective. The map $M_k\to \Theta_s(M_k)$ furthur induces an injective map $M_k\inj \tilde{M}_k$.
\begin{proof}
The proof is similar to the proof of Theorem \ref{thm:quo} and Corollary \ref{cor:quo}.
\end{proof}
\end{theorem}

Now we aim to show that the extension inside $\tilde{M}$ is non-split and give the universal extension.
\begin{lemma}\label{lem:nonsplitttt}
The exact sequence 
\begin{align}
    0\to M_{c,k}\to \tilde{M}_k\to M_{\alg,k}^{\oplus 2}\to 0
\end{align}
is non-split. Moreover, for any non-zero map $\alpha:\tilde{M}_k\to M_{\alg,k}$, the kernel $\ker\alpha$ is a non-split extension of $M_{\alg,k}$ by $M_{c,k}$.
\begin{proof}
Using the translation functor $T_{(0,0)}^{(k,0)}$ we may reduce to the case $k=0$. To prove the non-split result we need the geometric realization of $\tilde{M}_k$. Recall $\tilde{M}:=\tilde{M}_0=\Hom_{E[D_L^\times]}(\tau,H^0(\Sigma_n,\Omega_{\Sigma_n}^{1}))$ and $M_c:=M_{c,0}=\Hom_{E[D_L^\times]}(\tau,H^0(\Sigma_n,\calO_{\Sigma_n}))$. As $d:\calO_{\Sigma_n}\to \Omega_{\Sigma_n}^1$ is given by $f\mapsto u^+(f)dz$ up to scalars (see \cite[Lemme 3.5]{DLB17} for example), where $z$ is a local coordinate on $\Sigma_n$ such that $u^+(z)=1$, we see $u^+$ is injective on $\tilde{M}$. Suppose the extension 
\begin{align}
    0\to M_{c}\to \tilde{M}\to M_{\alg,0}^{\oplus 2}\to 0
\end{align}
splits $\GL_2(L)$-equivariantly, then $\tilde{M}$ contains a dual of a smooth representation $M_{\alg,0}$. This shows $u^+$ is not injective on $\tilde{M}$, which is a contradiction. Similarly one can show the kernel $\ker\alpha$ is a non-split extension of $M_{\alg,k}$ by $M_{c,k}$ for any non-zero map $\alpha:\tilde{M}_k\to M_{\alg,k}$.
\end{proof}
\end{lemma}

\begin{remark}
In fact, one can show the non-splitness inside $\tilde{M}_k$ using abstract non-sense. Consider $\Hom_{D_\sigma(G,E)}(M_{\alg,k},\tilde{M}_k)$. By Proposition \ref{prop:Hom}, we know the natural map 
\begin{align}
    \Hom_{D_\sigma(G,E)}(M_{\alg,k},\tilde{M}_k)\to \Hom_{D_\sigma(G,E)}(M_{\alg,k},\Theta_s(M_{c,k}))
\end{align}
is an isomorphism. But $\Hom_{D_\sigma(G,E)}(M_{\alg,k},\Theta_s(M_{c,k}))\isom  \Hom_{D_\sigma(G,E)}(T^{(-1,0)}_{(k,0)} M_{\alg,k},T^{(-1,0)}_{(k,0)}(M_{c,k}))=0$ as $T^{(-1,0)}_{(k,0)} M_{\alg,k}=0$.
\end{remark}



Now we can calculate the dimension of $\Ext^1_{D_\sigma(G,E)}(M_{\alg,k},M_{c,k})$. The way we calculate this first extension group is by using the wall-crossing method to relate it to certain Hom spaces. Let $\alpha:\tilde{M}_k\to M_{\alg,k}$ be a non-zero map. We know $\alpha$ is surjective, and by Proposition \ref{prop:Hom}, $\ker\alpha$ is an extension of $M_{\alg,k}$ by $M_{c,k}$. If $\alpha$ is zero, then we simply assign it to the zero element in $\Ext^1_{D_\sigma(G,E)}(M_{\alg,k},M_{c,k})$. Hence we get an $E$-linear map 
\begin{align}
    \Phi:\Hom_{D_\sigma(G,E)}(\tilde{M}_k, M_{\alg,k})\to \Ext^1_{D_\sigma(G,E)}(M_{\alg,k},M_{c,k}), \alpha\mapsto \ker\alpha.
\end{align}

\begin{theorem}\label{thm:Phiisom}
The map $\Phi$ is an isomorphism.
\begin{proof}
First we show $\Phi$ is surjective. This actually follows from Theorem \ref{thm:inj}. Indeed, let $M_k$ be a $D_\sigma(G,E)$-module, which is  non-split extension of $M_{\alg,k}$ by $M_{c,k}$. Then by wall-crossing we get an injective map $M_k\inj \tilde{M}_k$. Hence $M_k$ is the kernel of $\tilde{M}_k\to \tilde{M}_k/M_k\isom M_{\alg,k}$.

Then we show $\Phi$ is injective. First of all, Lemma \ref{lem:nonsplitttt} shows a non-zero map $\alpha:\tilde{M}_k\to M_{\alg,k}$ is mapped to a non-zero element in $\Ext^1_{D_\sigma(G,E)}(M_{\alg,k},M_{c,k})$. Let $\alpha,\alpha'$ be non-zero maps from $\tilde{M}_k$ to $M_{\alg,k}$ such that $\ker\alpha$ and $\ker\alpha'$ gives the same extension class. In this way we obtain a morphism of exact sequence with the first and the third vertical map being identity:
\begin{center}
\small
\begin{tikzpicture}[descr/.style={fill=white,inner sep=1.5pt}]
    \matrix (m) [
        matrix of math nodes,
        row sep=2.5em,
        column sep=2.5em,
        text height=1.5ex, 
        text depth=0.25ex
    ]
    { 0 & M_{c,k} & M_k & M_{\alg,k} & 0  \\
    0 & M_{c,k} & M_k' & M_{\alg,k} & 0  \\
    };
    \path[->,font=\scriptsize]
    (m-1-1) edge (m-1-2)
    (m-1-2) edge (m-1-3)
    (m-1-3) edge (m-1-4)
    (m-1-4) edge (m-1-5)
    (m-2-1) edge (m-2-2)
    (m-2-2) edge (m-2-3)
    (m-2-3) edge (m-2-4)
    (m-2-4) edge (m-2-5)
    ;
    \path[->,font=\scriptsize]
    (m-1-3) edge node[right] {$\sim$} node[left] {$\phi$} (m-2-3)
    ;
    \draw[double distance = 1.5pt] 
    (m-1-2) -- (m-2-2)
    ;
    \draw[double distance = 1.5pt] 
    (m-1-4) -- (m-2-4)
    ;
\end{tikzpicture}. 
\end{center}
As $M_k$, $M_k'$ are non-split extensions by Lemma \ref{lem:nonsplitttt}, we get induced maps $\beta:M_k\inj \tilde{M}$ and $\beta':M_k'\inj \tilde{M}$. As the isomorphism $\phi:M_k\aisom M_k'$ is identity on $M_{c,k}$, from the construction we see $\phi$ induces the identity map on $\tilde{M}\surj \tilde{M}/\im \beta$ and $\tilde{M}\surj \tilde{M}/\im \beta'$. This shows that $\Phi$ is injective, and finishes the proof.
\end{proof}
\end{theorem}

\begin{corollary}\label{cor:Ext1onlyonDr}
We have $\dim_E\Ext^1_{D_\sigma(G,E)}(M_{\alg,k},M_{c,k})=2$.
\begin{proof}
By Theorem \ref{thm:Phiisom}, we know 
\begin{align}
    \dim_E\Ext^1_{D_\sigma(G,E)}(M_{\alg,k},M_{c,k})=\dim_E\Hom_{D_\sigma(G,E)}(\tilde{M}_k, M_{\alg,k}).
\end{align}
As $\Hom_{D_\sigma(G,E)}(M_{c,k},M_{\alg,k})=0$ by Proposition \ref{prop:Hom}, we see that $\dim_E\Hom_{D_\sigma(G,E)}(\tilde{M}_k, M_{\alg,k})=2$.
\end{proof}
\end{corollary}

\begin{remark}\label{rmk:Phiisom2}
Using the same method, one can similarly show that $\dim_E\Ext^1_{D_\sigma(G,E)}(M_{c,k},M_{\alg,k})=2$.
\end{remark}

\subsection{Extension classes and isomorphism classes}
Let $M$ be a $D_\sigma(G,E)$-module which is an extension of $M_{\alg,k}$ by $M_{c,k}$. We want to show that the isomorphism class of $M$ determines the extension class of $M$ inside $\Ext^1_{D_\sigma(G,E)}(M_{\alg,k},M_{c,k})$ in up to scalars. A key step is the following result.
\begin{proposition}\label{prop:END}
We have $\End_{D_\sigma(G,E)}(M_{c,0})=E$.
\begin{proof}
Recall that $\Omega$ is the Drinfeld upper half plane, viewed as an adic space over $E$ such that $L$ acts via $\sigma$. Let $\calO(\Omega)$ be the global section of the structure sheaf on $\Omega$. As the image of $\pi_{\GM,n}:\Sigma_n\to \fl$ is $\Omega$, the algebra $\calO(\Omega)$ acts on $\calO(\Sigma_n)$. As $D_L^\times$ acts trivially on $\Omega$, we see $M_{c,0}$ also carries an action of $\calO(\Omega)$. First we want to show that every continuous $\GL_2(L)$-equivariant endomorphism on $M_{c,0}$ is $\calO(\Omega)$-semi-linear. Let $\Phi$ be a continuous $\GL_2(L)$-equivariant endomorphism on $M_{c,0}$. From the proof of \cite[Proposition 9.16]{DLB17}, we see that it suffices to prove that $\Phi$ commutes with multiplication by $z\in\calO(\Omega)$, where $z$ is a local coordinate on $\Omega$.  Let $a^+=\left( \begin{matrix} 1&0\\0&0\end{matrix} \right)\in\gl_2(L)\ox_{L,\sigma}E$. By \cite[Lemme 3.5]{DLB17} we know that $a^+-1=u^+z$ on $\calO(\Sigma_n)$. As $\Phi$ is $\GL_2(L)$-equivariant, $\Phi$ commutes with the action of $a^+$ and $u^+$. From this we see that $u^+\Phi (zs)=u^+z\Phi(s)$ where $s\in M_{c,0}$. As $d:\calO_{\Sigma_n}\to \Omega^1_{\Sigma_n}$ is (up to scalar) given by $f\mapsto u^+(f)dz$ and the map $M_{c,0}\to \tilde{M}_0$ induced by $d$ is injective, we see $u^+$ is injective on $M_{c,0}$ and hence $\Phi (zs)=z\Phi(s)$, which shows $z$ commutes with $\Phi$. Then the continuity of $\Phi$ implies that 
\begin{align}
    \End_{D_\sigma(G,E)}(M_{c,0})=\End_{D_\sigma(G,E),\calO(\Omega)}(M_{c,0}).
\end{align}
Here, $\End_{D_\sigma(G,E),\calO(\Omega)}(-)$ means continuous endomorphisms linear with respect to the action of the twisted group algebra $\calO(\Omega)[G]$. Next we show $\End_{D_\sigma(G,E),\calO(\Omega)}(M_{c,0})=E$.
We use the results in \cite{Koh11LT} on Lubin-Tate bundles and Drinfeld bundles. See \cite{Koh11LT} for definitions of Lubin-Tate bundles and Drinfeld bundles and the functor $D_{\LT}(-)$. As the coherent sheaf associated to $M_{c,0}=\Hom_{D_L^\times}(\tau,\calO_{\Sigma_n}(\Sigma_n))$ is a Lubin-Tate bundle by \cite[Theorem 4.5]{Koh11LT}, and $D_{\LT}(-)$ is fully faithful on Lubin-Tate bundles, we deduce
\begin{align}
    \Hom_{D_\sigma(\GL_2(L),E),\calO(\Omega)}(M_{c,0},M_{c,0})\isom \Hom_{\calO_{\check{\fl}},D^\times_L}(\calO_{\check{\fl}}\ox_E \tau,\calO_{\check{\fl}}\ox_E \tau)\isom H^0(\check{\fl},\calO_{\check{\fl}})\ox_E \End_{D^\times_L}(\tau)\isom E.
\end{align}
where $\check{\fl}\isom \bbP^1$ is the period domain for the Gross-Hopkins period map on the Lubin-Tate curves.
\end{proof}
\end{proposition}


Recall that from Proposition \ref{prop:basicoftranslation} we can deduce $M_{c,k}\isom T_{(0,0)}^{(k,0)}M_{c,0}$ for $k\ge 0$. We extend this definition to $k=-1$ by letting $M_{c,-1}:=T_{(0,0)}^{(-1,0)}M_{c,0}$.
\begin{corollary}\label{cor:END}
For any $k\ge -1$, we have $\End_{D_\sigma(G,E)}(M_{c,k})=E$.
\begin{proof}
Let $k\ge 0$ be an integer. As $M_{c,k}\isom T_{(0,0)}^{(k,0)}M_{c,0}$ and $T_{(0,0)}^{(k,0)}$ is biadjoint to $T_{(k,0)}^{(0,0)}$, we have 
\begin{align}
    \Hom_{D_\sigma(G,E)}(M_{c,k},M_{c,k})=\Hom_{D_\sigma(G,E)}(M_{c,0},T_{(k,0)}^{(0,0)}T_{(0,0)}^{(k,0)}M_{c,0})=\Hom_{D_\sigma(G,E)}(M_{c,0},M_{c,0})=E
\end{align}
as $T_{(k,0)}^{(0,0)}T_{(0,0)}^{(k,0)}M_{c,0}\isom M_{c,0}$. When $k=-1$, we have 
\begin{align}
    \Hom_{D_\sigma(G,E)}(T_{(0,0)}^{(-1,0)}M_{c,0},T_{(0,0)}^{(-1,0)}M_{c,0})=\Hom_{D_\sigma(G,E)}(M_{c,0},\Theta_s M_{c,0})
\end{align}
As $M_{c,0}$ has trivial infinitesimal character, we have 
\begin{align}
    \Hom_{D_\sigma(G,E)}(M_{c,0},\Theta_s M_{c,0})=\Hom_{D_\sigma(G,E)}(M_{c,0},\Theta_s M_{c,0}[\calZ(\frg_\sigma)=\tilde\chi_{(0,0)}])
\end{align}
By Lemma \ref{lem:wcsemisimple} and Proposition \ref{prop:Hom}, we deduce 
\begin{align}
    \Hom_{D_\sigma(G,E)}(M_{c,0},\Theta_s M_{c,0}[\calZ(\frg_\sigma)=\tilde\chi_{(0,0)}])=\Hom_{D_\sigma(G,E)}(M_{c,0},M_{c,0}-M_{\alg,0}^{\oplus 2})=\Hom_{D_\sigma(G,E)}(M_{c,0},M_{c,0})=E
\end{align}
where the last step follows from Proposition \ref{prop:END}.
\end{proof}
\end{corollary}

\begin{proposition}\label{prop:isompext}
Let $M$, $M'$ be two non-split extensions of $M_{\alg,k}$ by $M_{c,k}$, which are isomorphic to each other as $D_\sigma(G,E)$-modules. Then $M$, $M'$ have the same image in $\bbP\Ext^1_{D_\sigma(G,E)}(M_{\alg,k},M_{c,k})$. Here $\bbP\Ext^1_{D_\sigma(G,E)}(M_{\alg,k},M_{c,k})$ is the quotient $(\Ext^1_{D_\sigma(G,E)}(M_{\alg,k},M_{c,k})\bs\{0\})/E^\times$, where $0$ stands for the trivial extension class.
\begin{proof}
Let $M$, $M'$ be two extensions of $M_{\alg,k}$ by $M_{c,k}$, such that $M\isom M'$ are isomorphic as $D_\sigma(G,E)$-modules. Suppose $\phi:M\aisom M'$ is a $D_\sigma(G,E)$-equivariant isomorphism. By considering the wall-crossing of $\phi$, we get a morphism of exact sequence 
\begin{center}
\small
\begin{tikzpicture}[descr/.style={fill=white,inner sep=1.5pt}]
    \matrix (m) [
        matrix of math nodes,
        row sep=2.5em,
        column sep=2.5em,
        text height=1.5ex, 
        text depth=0.25ex
    ]
    { 0 & M & \tilde M & M_{\alg,k} & 0  \\
      0 & M' & \tilde M & M_{\alg,k} & 0 \\
    };
    \path[->,font=\scriptsize]
    (m-1-1) edge (m-1-2)
    (m-1-2) edge (m-1-3)
    (m-1-3) edge (m-1-4)
    (m-1-4) edge (m-1-5)
    (m-2-1) edge (m-2-2)
    (m-2-2) edge (m-2-3)
    (m-2-3) edge (m-2-4)
    (m-2-4) edge (m-2-5)
    ;
    \path[->,font=\scriptsize]
    (m-1-4) edge node[right] {$\sim$} (m-2-4)
    (m-1-2) edge node[right] {$\phi$} (m-2-2)
    (m-1-3) edge node[right] {$\tilde\phi$} (m-2-3)
    ;
\end{tikzpicture}
\end{center}
where the middle arrow comes from a $D_\sigma(G,E)$-equivariant automorphism of $T_{(k,0)}^{(-1,0)}M_{c,k}$. By Corollary \ref{cor:END}, we know that $\End_{D_\sigma(G,E)}(M_{c,-1})=E$, so that $T_{(0,0)}^{(-1,0)}\phi$ is given by a scalar. From this we see that the middle isomorphism $\tilde\phi$ is given by multiplying a scalar in $E^\times$. Using this, one can construct an isomorphism of extension classes of $\pi$ and $\pi'$ up to multiplying by scalars as follows. After multiplying by a scalar on the isomorphism $\phi:M\aisom M'$, we can ensure this isomorphism induces the identity map on $M_{\alg,k}$. Besides, as $\End_{D_\sigma(G,E)}(M_{c,k})=E$, we can modify the inclusion map $M_{c,k}\inj M$ so that $\phi|_{M_{c,k}}$ is also given by the identity map. Hence $M$ and $M'$ are also isomorphic as extensions of $M_{\alg,k}$ by $M_{c,k}$ up to scalars.
\end{proof}
\end{proposition}

\section{Applications to the $p$-adic local Langlands program for $\GL_2(L)$}

\subsection{Generalities on Deligne-Fontaine modules and Hodge filtrations}
In this section, we recall some generalities on Deligne-Fontaine modules, Weil-Deligne representations and Hodge filtrations. The key point is to track how these data interact with respect to the decomposition $L\ox_{\bbQ_p}E\isom \prod_{\sigma\in \Sigma}E$. Our main reference is \cite[2.2]{Bre19}.

Let $L'$ be a finite Galois extension of $L$ and $L'_0$ its maximal unramified subextension of $\bbQ_p$. We assume $E$ is sufficiently large in the sense that $[L_0':\bbQ_p]=|\Hom_{\bbQ_p}(L_0',E)|$. Consider the following two categories:
\begin{enumerate}[(i)]
    \item Let $\WD_{L'/L}$ be the category of Weil-Deligne representations $(r,N,V)$, where $V$ is a finite dimensional $E$-vector space, $r$ is a smooth $E$-linear action of ${{\text{W}_L}}$ (the Weil group of $L$) on $V$ which is unramified when restricted to $W_{L'}$, $N$ is an $E$-linear endomorphism of $V$, such that $r(g)Nr(g)^{-1}=||g||N$ where $||\cdot||:{{\text{W}_L}}\to q^{\bbZ}$ is the unramified character such that $||\Frob||=q=p^f$ with $\Frob$ the geometric Frobenius.
    \item Let $\DF_{L'/L}$ be the category of Deligne-Fontaine modules $(D,\vphi,N,\Gal(L'/L))$ where $D$ is a free $L_0'\ox_{\bbQ_p}E$-module of finite rank endowed with a $L_0'$-semilinear, $E$-linear Frobenius $\vphi:D\to D$, and $L_0'\ox_{\bbQ_p}E$-linear endomorphism $N:D\to D$ such that $N\vphi=p\vphi N$ and $L_0'$-semilinear, $E$-linear action of $\Gal(L'/L)$ on $D$ commuting with $\vphi$ and $N$. 
\end{enumerate} 
These two categories are known to be equivalent, compatible with $L'$. More precisely, let $\ul D\in \text{DF}_{L'/L}$ be a Deligne-Fontaine module, and $\sigma_0':L_0'\inj E$ be an embedding. Then it is known that $D_{L_0',\sigma_0'}:=D\ox_{L_0'\ox_{\bbQ_p}E,\sigma_0'\ox\id}E$ can be equipped with an action of the Weil-Deligne group.

Let $D_{L'}:=D\ox_{L_0'}L'$ and we extend the action of $\Gal(L'/L)$ by $g((l'\ox e).d)=(g(l')\ox e).g(d)$ with $l'\in L'$, $e\in E$, and $d\in D$.  The decomposition $L\ox_{\bbQ_p}E\aisom\prod_{\sigma\in \Sigma} E$, $l\ox e\mapsto (\sigma(l)e)_{\sigma}$ induces a decomposition $D_{L'}\isom \bigoplus_{\sigma\in \Sigma} D_{L',\sigma}$, where each $D_{L',\sigma}$ is a $L'\ox_{L,\sigma}E$-module free with same rank as $D$ over $L_0'\ox_{\bbQ_p}E$, and the action of $\Gal(L'/L)$ on $D_{L'}$ induces an action of $D_{L',\sigma}$. By Galois descent, we get an $L'\ox_{L,\sigma}E$-linear isomorphism for each $\sigma$:
\begin{align}
    L'\ox_L D_{L',\sigma}^{\Gal(L'/L)}\aisom D_{L',\sigma}.
\end{align}
Let $\sigma':L'\inj E$ be an embedding such that $\sigma'|_{L}=\sigma$, and define $D_{L',\sigma'}=D_{L'}\ox_{L'\ox_{\bbQ_p}E,\sigma'\ox\id}E$. The injection $D\inj D_{L'}$ induces an isomorphism of $E$-vector spaces $D_{L_0',\sigma'|_{L_0'}}\aisom D_{L',\sigma'}$, so that we get a Weil-Deligne action on $D_{L',\sigma'}$ by transporting structures. Moreover, the composition 
\begin{align}
    D_{L',\sigma}^{\Gal(L'/L)}\inj D_{L',\sigma}\surj D_{L',\sigma'}
\end{align}
is an isomorphism of $E$-vector spaces, so that by transporting structure we also get a Weil-Deligne action on $D_{L',\sigma}^{\Gal(L'/L)}$. By \cite[Lemme 2.2.1]{Bre19}, we see that up to conjugation, the Weil-Deligne action on $D_{L',\sigma}^{\Gal(L'/L)}$ is independent of the choice of $\sigma':L'\inj E$ such that $\sigma'|_L=\sigma$. From now on we fix such an embedding $\sigma'$.


Now suppose $V$ is an $n$-dimensional $E$-linear de Rham representation of $\Gal(\bar L/L)$. If $L'/L$ is a finite Galois extension such that $V|_{\Gal(\bar L/L')}$ is semi-stable, then we can associate the Deligne-Fontaine module $\ul D=((B_{\st}\ox_{\bbQ_p}V)^{\Gal(\bar L/L')},\vphi,N,\Gal(L'/L))$, with the Hodge filtration
\begin{align}
    D_{L'}^{\Gal(L'/L)}=(B_{\dR}\ox_{\bbQ_p}V)^{\Gal(\bar L/L)}\supset \Fil^i D_{L'}^{\Gal(L'/L)}=(\Fil^iB_{\dR}\ox_{\bbQ_p}V)^{\Gal(\bar L/L)}).
\end{align}
Via the decomposition $D_{L'}=\bigoplus_{\sigma\in\Sigma}D_{L',\sigma}$, we get the $\sigma$-Hodge filtration
\begin{align}
    (\Fil^i B_{\dR}\ox_{\bbQ_p}V)^{\Gal(\bar L/L)}\ox_{L\ox_{\bbQ_p}E,\sigma\ox \id}E\subset D_{L',\sigma}^{\Gal(L'/L)}
\end{align}
On the Weil-Deligne representation $D_{L',\sigma}^{\Gal(L'/L)}$ for each $\sigma\in \Sigma$.


Conversely, let $\ul D=(D,\vphi,N,\Gal(L'/L))\in \text{DF}_{L'/L}$ be a Deligne-Fontaine module of rank $n$ over $L_0'\ox_{\bbQ_p}E$, with an admissible filtration $(D_{L'},\Fil^{\bullet}D_{L'})$. For example see \cite[D\'efinition 7.1]{Bre16Soc1} for definition of admissiblilty of the filtration. Set $V_{\st}^0(\ul D,\Fil^{\bullet})=(B_{\st}\ox_{L_0'}D)^{\vphi=1,N=0}$ and $V_{\st}^1(\ul D,\Fil^{\bullet})=B_{\dR}\ox_{L'}D_{L'}/\Fil^0(B_{\dR}\ox_{L'}D_{L'})$. After fixing a branch of the $p$-adic logarithm, the inclusion $B_{\st}\inj B_{\dR}$ induces a map 
\begin{align}
    V^0_{\st}(\ul D,\Fil^{\bullet})\to V^1_{\st}(\ul D,\Fil^{\bullet})
\end{align}
which is $\Gal(\bar L/L)$-equivariant for the diagonal action. By the main result of \cite{CF00}, this map is surjective with kernel $V_{\st}(\ul D,\Fil^{\bullet})=\Fil^0(B_{\st}\ox_{L_0'}D)^{\vphi=1,N=0}$, which is an $n$-dimensional $E$-linear continuous representation of $\Gal(\bar L/L)$. This procedure induces an equivalence of categories between $n$-dimensional potentially semi-stable representations of $\Gal(\bar L/L)$ over $E$ such that stable via restriction to $\Gal(\bar L/L')$ and Deligne-Fontaine modules of rank $n$ over $L_0'\ox_{\bbQ_p}E$ with admissible Hodge filtrations.

The above construction extracts information from potentially semi-stable representations of $\Gal(\bar L/L)$ and encodes them into linear-algebraic structures. Compared to Weil-Deligne representations, which are among the main objects considered in the classical local Langlands correspondence, potentially semi-stable representations contain additional information, namely, the Hodge filtrations. Furthermore, although there is an equivalence of categories between Weil-Deligne representations and Deligne-Fontaine modules, strictly speaking, Weil-Deligne representations only use the information from $\vphi^f$ with $f=[L_0:\bbQ_p]$. Through the isomorphisms $\vphi^d: D_{L_0', \sigma_0'} \xrightarrow{\sim} D_{L_0', \sigma_0' \circ \vphi_0'^{-d}}$ (up to an $E$-linear automorphism of $D_{L_0', \sigma_0'}$ which is independent of Hodge filtrations) of Weil-Deligne representations for $d \in \mathbb{Z}$, and $D_{L_0', \sigma'|_{L_0'}} \xrightarrow{\sim} D_{L', \sigma}^{\Gal(L'/L)}$ for each $\sigma \in \Sigma$, the $\sigma$-Hodge filtration on $D_{L', \sigma}^{\Gal(L'/L)}$ also defines a filtration on $D_{L', \tau}^{\Gal(L'/L)}$ for each $\tau \neq \sigma$. When $\End_{\text{DF}}(\underline{D}) = E$, these new filtrations are completely determined by the original filtrations. In this case, we see the Galois representation $V$ is completely determined by the associated Weil-Deligne representation, and $\sigma$-Hodge filtrations for each $\sigma \in \Sigma$.

\subsection{Pro-\'etale cohomology of Drinfeld towers of dimension $1$}\label{sec:prelim}
Finally, we apply our results to study the pro-\'etale cohomology of Drinfeld towers of dimension $1$. For the convenience of the reader and consistency of notation, we (re)define some important notations. Recall that $\breve{\calM}_{\Dr,n}$ is the Drinfeld tower of dimension $1$ over $\breve{L}$ of level $n$. We base change them to $C$ and we denote them by $\calM_{\Dr,n}:=\breve{\calM}_{\Dr,n}\times_{\breve{L}}C$.  

Following notation in \cite{CDN20}, let $M$ be a Deligne-Fontaine module over $E\ox_{\bbQ_p}L_0^{\ur}$, which is a $(\vphi,N,\Gal(\bar L/L))$-module, such that $\vphi$ is $E$-linear, $L_0$-Frobenius semi-linear, $N$ is $E\ox_{\bbQ_p}L_0^{\ur}$-linear and $N\vphi=p\vphi N$, and the action of $\Gal(\bar L/L)$ is smooth and $L_0$-semilinear. We assume that $M$ is of slope $\frac{1}{2}$. To $M$, we can associate following representations:
\begin{itemize}
    \item Let $\WD(M)$ be the Weil-Deligne representation corresponding to $M$. This is a $2$-dimensional $E$-linear representation of $\WD_L$, the Weil-Deligne group of $L$. We assume $\WD(M)$ is irreducible.
    \item Let $\LL(M)$ be the smooth irreducible representation of $\GL_2(L)$ corresponding to $\WD(M)$ via the local Langlands correspondence. As $\WD(M)$ is irreducible, $\LL(M)$ is supercuspidal.
    \item Let $\JL(M)$ be the smooth irreducible representation of $D_L^\times$ corresponding to $\LL(M)$ via the local Jacquet-Langlands correspondence. As $\LL(M)$ is supercuspidal, the action of $D_L^\times$ on $\JL(M)$ does not factor through the reduced norm.  After twisting by characters, we may assume that $\varpi$ acts trivially on $\JL(M)$. From the discussion in \cite[5.1]{CDN20}, we see that this assumption does not affect the main results.
\end{itemize}
We also need some spaces to construct $p$-adic Galois representations from $M$.
\begin{itemize}
    \item Let $M_{\dR}:=(M\ox_{L_0^{\ur}}C)^{\Gal(\bar L/L)}$, which is a free $E\ox_{\bbQ_p}L$-module of rank $2$.
    \item Define $X_{\st}^+(M)=(B_{\cris}^+\ox_{\bbQ_p^{\ur}}M)^{\vphi=p}=(B_{\cris}^+\ox_{\bbQ_p^{\ur}}M[-1])^{\vphi=1}$, which is a Banach-Colmez Space. Here the twist $[m]$ for $m\in\bbZ$ denotes multiplication of the $\vphi$-action by $p^m$.
\end{itemize}
Let $\calL\subset M_{\dR}$ be a locally free $L\ox_{\bbQ_p}E$-module of rank $1$. Define 
\begin{align}
    V_{M,\calL}:=\ker(X_{\st}^+(M)\to B_{\dR}^+\ox_L M_{\dR}\ov{\theta}\to  C\ox_L (M_{\dR}/\calL)).
\end{align}
By \cite{CF00}, we know $V_{M,\calL}$ is a $2$-dimensional $E$-linear de Rham representation of parallel Hodge-Tate weight $0,1$, such that $D_{\pst}(V_{M,\calL})\isom M[-1]$.  Again, we note that our normalization is such that the Hodge-Tate weight of the $p$-adic cyclotomic character is $1$. Conversely, suppose that $V$ is a $2$-dimensional $E$-linear de Rham representation of parallel Hodge-Tate weight $0,1$ such that $D_{\pst}(V_{M,\calL})\isom M[-1]$, then by the admissibility condition $V\isom V_{M,\calL}$ for some unique rank 1 locally free submodule $\calL\subset M_{\dR}$.

Let $H^1_{\dR}(\calM_{\Dr,n})$ be the first de Rham cohomology of $\calM_{\Dr,n}$, we define 
\begin{align}
    H^1_{\dR}(\calM_{\Dr,\infty}):=\dlim_n H^1_{\dR}(\calM_{\Dr,n}).
\end{align}
This is an ind-Frech\'et space over $C$. We also have a $E\ox_{\bbQ_p}C$-linear $\GL_2(L)$-equivariant isomorphism 
\begin{align}
    \Hom_{E[D_L^\times]}(\JL(M), E\ox_{\bbQ_p}H^1_{\dR}(\calM_{\Dr,\infty}))=C\hat\ox_L M_{\dR}\hat\ox_E \LL(M)^*.
\end{align}
Here, $\LL(M)^*$ denotes the strong dual of $\LL(M)$.

Let $H^1_{\HK}(\calM_{\Dr,n})$ be the Hyodo-Kato cohomology of $\calM_{\Dr,n}$, and we define 
\begin{align}
    H^1_{\HK}(\calM_{\Dr,\infty}):=\dlim_n H^1_{\HK}(\calM_{\Dr,n}).
\end{align}
This is an ind-Frechet space over $\breve L_0$, equipped with semi-linear $\vphi$-action, an operator $N$ such that $N\vphi=p\vphi N$, and a semi-linear smooth $\Gal(\bar L/L)$-action. This space also carries a $\GL_2(L)\times D_L^\times$-action. Then we have an isomorphism of $(\vphi,N,\Gal(\bar L/L))$-modules over $E\ox_{\bbQ_p}\breve L_0$:
\begin{align}
    \Hom_{E[D_L^\times]}(\JL(M),E\ox_{\bbQ_p}H^1_{\HK}(\calM_{\Dr,\infty}))\aisom \breve{L}_0\ox_{L_0^{\ur}}M\hat\ox_E\LL(M)^*
\end{align}
which is also compatible with the $\GL_2(L)$-action. Besides, we have a natural isomorphism of ind-Frechet spaces over $C$:
\begin{align}
    \iota_{\HK}:C\hat\ox_{\breve L_{0}}H^1_{\HK}(\calM_{\Dr,\infty})\aisom H^1_{\dR}(\calM_{\Dr,\infty})
\end{align}
which is equivariant for the $\GL_2(L)\times D_L^\times$-action. This isomorphism induces a commutative diagram (\cite[Th\'eor\`eme 0.4]{CDN20}) 
\begin{center}
\begin{tikzpicture}[descr/.style={fill=white,inner sep=1.5pt}]
    \matrix (m) [
        matrix of math nodes,
        row sep=2.5em,
        column sep=2.5em,
        text height=1.5ex, 
        text depth=0.25ex
    ]
    { \Hom_{E[D_L^\times]}(\JL(M),E\ox_{\bbQ_p}H^1_{\HK}(\calM_{\Dr,\infty}))&M\ox_E\LL(M)^*  \\
    \Hom_{E[D_L^\times]}(\JL(M), E\ox_{\bbQ_p}H^1_{\dR}(\calM_{\Dr,\infty}))&C\hat\ox_L M_{\dR}\hat\ox_E \LL(M)^* \\
    };
    \path[->,font=\scriptsize]
    (m-1-1) edge node[auto] {$\sim$} (m-1-2)
    (m-2-1) edge node[auto] {$\sim$} (m-2-2)
    ;
    \path[->,font=\scriptsize]
    (m-1-1) edge node[auto] {$\iota_{\HK}$} (m-2-1)
    (m-1-2) edge node[auto] {} (m-2-2)
    ;
\end{tikzpicture}
\end{center}
where the left vertical map is induced by the Hyodo-Kato map, and the right vertical map is induced by the natural embedding $M\subset M_{\dR}\ox_{L}C$. The vertical maps are isomorphisms if we tensor the first row with $-\hat\ox_{\breve L_0}C$. The diagram is $\GL_2(L)$-equivariant, and the first horizontal map is also compatible with the $(\vphi,N,\Gal(\bar L/L))$-action. We decompose the above diagram with respect to the decompositon 
\begin{align}
    E\ox_{\bbQ_p}C\isom\prod_{\sigma\in\Sigma}E\ox_{\sigma,L,\iota}C
\end{align}
and pick the $\sigma$-component, so that to get a commutative diagram
\begin{center}
\begin{tikzpicture}[descr/.style={fill=white,inner sep=1.5pt}]
    \matrix (m) [
        matrix of math nodes,
        row sep=2.5em,
        column sep=2.5em,
        text height=1.5ex, 
        text depth=0.25ex
    ]
    { \Hom_{E[D_L^\times]}(\JL(M),E\ox_{\sigma,L,\iota}C\hat\ox_{\breve L_0}H^1_{\HK}(\calM_{\Dr,\infty}))&C\hat \ox_{\iota,L,\sigma}\WD(M)\ox_E\LL(M)^*  \\
    \Hom_{E[D_L^\times]}(\JL(M), E\ox_{\sigma,L,\iota}H^1_{\dR}(\calM_{\Dr,\infty}))&C\hat\ox_{\iota,L,\sigma} M_{\dR,\sigma}\hat\ox_E \LL(M)^* \\
    };
    \path[->,font=\scriptsize]
    (m-1-1) edge node[auto] {$\sim$} (m-1-2)
    (m-2-1) edge node[auto] {$\sim$} (m-2-2)
    ;
    \path[->,font=\scriptsize]
    (m-1-1) edge node[auto] {$\iota_{\HK}$} (m-2-1)
    (m-1-2) edge node[auto] {$\isom$} (m-2-2)
    ;
\end{tikzpicture}.
\end{center}
Here, $M_{\dR,\sigma}=M_{\dR}\ox_{L\ox_{\bbQ_p}E}L\ox_{L,\sigma}E$. Indeed, there exists a sufficiently large finite extension $L'$ of $L$, such that $\Gal(\bar L/L')$ acts trivially on $M$. Let $L_0'$ be the maximal unramified subextension of $\bbQ_p$ inside $L'$, then there exists a $(\vphi,N,\Gal(L'/L))$-module $M_{L_0'}$ over $E\ox_{\bbQ_p}L_0'$, such that $M\isom M_{L_0'}\ox_{L_0'}L_0^{\ur}$. Let $M_{L'}\isom M\ox_{L_0'}L'$, and we pick an embedding $L'\inj C$ extending the embedding $\iota:L\inj C$. Then we have 
\begin{align}
    C\ox_{L_0^{\ur}}M\isom C\ox_{L_0'}M_{L_0'}\isom C\ox_{L'} M_{L'}\isom C\ox_{L}M_{L'}^{\Gal( L'/L)}
\end{align}
and hence 
\begin{align}
    C\ox_{L_0^{\ur}}M\ox_{C\ox_{\bbQ_p}E}C\ox_{L,\sigma}E\isom C\ox_{L}M_{L',\sigma}^{\Gal( L'/L)}
\end{align}
where $M_{L',\sigma}^{\Gal( L'/L)}=M_{L'}^{\Gal( L'/L)}\ox_{E\ox_{\bbQ_p}L}E\ox_{L,\sigma}L$. From our discuss in the previous section, we know that $M_{L',\sigma}^{\Gal( L'/L)}$ can be equipped with an action of the Weil-Deligne group of $L$, and it is the associated Weil-Deligne representation of $M$ in \cite[2.2]{Bre19}.

The pro-\'etale cohomology of Drinfeld tower of dimension $1$ combines all above data together. Define $H^0(\calO_{\Dr,\infty}):=\dlim_n H^0(\calM_{\Dr,n},\calO_{\calM_{\Dr,n}})$ and $H^0(\Omega^1_{\Dr,\infty}):=\dlim_n H^0(\calM_{\Dr,n},\Omega^1_{\calM_{\Dr,n}})$. They are ind-Frech\'et spaces over $C$, carrying $\GL_2(L)\times D_L^\times$-action. Define $H^1_{\pro\et}(\calM_{\Dr,\infty},E(1)):=\dlim_nH^1_{\pro\et}(\calM_{\Dr,n},E(1))$, which carries $\GL_2(L)\times D_L^\times$-action and the $D_L^\times$-action is smooth. Here we note that in order to be consistent with the notation in \cite{CDN20}, we consider representations in $H^1_{\pro\et}(\calM_{\Dr,\infty},E(1))$ rather than $H^1_{\pro\et}(\calM_{\Dr,\infty},E)$. As $\calM_{\Dr,\infty}$ is over $C$, these two Galois representations differ only by a Tate twist:
\begin{align}
    H^1_{\pro\et}(\calM_{\Dr,\infty},E(1))\isom H^1_{\pro\et}(\calM_{\Dr,\infty},E)(1).
\end{align}
By \cite[Théorème 0.8]{CDN20}, the pro-\'etale cohomology of Drinfeld tower of dimension $1$ fits into the following diagram 
\begin{center}
\scriptsize
\begin{tikzpicture}[descr/.style={fill=white,inner sep=1.5pt}]
    \matrix (m) [
        matrix of math nodes,
        row sep=2.5em,
        column sep=2.5em,
        text height=1.5ex, 
        text depth=0.25ex
    ]
    { 0 & \Hom_{E[D_L^\times]}(\JL(M),E\ox_{\bbQ_p}H^0(\calO_{\Dr,\infty})) & \Hom_{E[D_L^\times]}(\JL(M),H^1_{\pro\et}(\calM_{\Dr,\infty},E(1))) & X_{\st}^+(M)\hat\ox_E \LL(M)^* & 0  \\
      0 & \Hom_{E[D_L^\times]}(\JL(M),E\ox_{\bbQ_p}H^0(\calO_{\Dr,\infty})) & \Hom_{E[D_L^\times]}(\JL(M),E\ox_{\bbQ_p}H^0(\Omega^1_{\Dr,\infty})) & C\ox_L M_{\dR}\hat\ox_E\LL(M)^* & 0 \\
    };
    \path[->,font=\scriptsize]
    (m-1-1) edge (m-1-2)
    (m-1-2) edge (m-1-3)
    (m-1-3) edge (m-1-4)
    (m-1-4) edge (m-1-5)
    (m-2-1) edge (m-2-2)
    (m-2-2) edge (m-2-3)
    (m-2-3) edge (m-2-4)
    (m-2-4) edge (m-2-5)
    ;
    \path[->,font=\scriptsize]
    (m-1-3) edge node[right] {} (m-2-3)
    (m-1-4) edge node[right] {$\theta$} (m-2-4)
    ;
    \draw[double distance = 1.5pt] 
    (m-1-2) -- (m-2-2)
    ;
\end{tikzpicture}.
\end{center}
The horizontal rows are exact, and the vertical rows are closed embeddings. The whole diagram is $E\ox_{\bbQ_p}C$-linear, $\GL_2(L)$-equivariant and $\Gal(\bar L/L)$-equivariant. The map $\theta$ in the above diagram is induced from $B_{\cris}^+\to B_{\dR}^+\to C$ and the natural isomorphism $C\ox_L M_{\dR}\isom C\ox_{L_0^{\ur}}M$.

Set 
\begin{align}
    \calO[M]&:=\Hom_{E[D_L^\times]}(\JL(M),E\ox_{\bbQ_p}H^0(\calO_{\Dr,\infty}))^{\Gal(\bar L/L)},\\
    \Omega^1[M]&:=\Hom_{E[D_L^\times]}(\JL(M),E\ox_{\bbQ_p}H^0(\Omega^1_{\Dr,\infty}))^{\Gal(\bar L/L)},\\
    H^1_{\pro\et}[M]&:=\Hom_{E[D_L^\times]}(\JL(M),H^1_{\pro\et}(\calM_{\Dr,\infty},E(1))).
\end{align}
Starting from the $E\ox_{\bbQ_p}L$-linear, $\GL_2(L)$-equivariant exact sequence
\begin{align}\label{eq:dRM}
    0\to \calO[M]\to \Omega^1[M]\to M_{\dR}\ox_E \LL(M)^*\to 0,
\end{align}
let $W_{M,\calL}'$ (this is not the dual of $W_{M,\calL}$) be the inverse image of $\calL\ox_E\LL(M)^*\subset M_{\dR}\ox_E \LL(M)^*$ inside $\Omega^1[M]$. Let $\calL^{-1}:=M_{\dR}^*/(M_{\dR}/\calL)^*$ such that under the natural bilinear map $M_{\dR}^*\ox M_{\dR}\to E\ox_{\bbQ_p}L$, it satisfies $\calL^{-1}\ox_{E\ox_{\bbQ_p}L}\calL\isom E\ox_{\bbQ_p}L$. Define $W_{M,\calL}$ to be the inverse image of $E\ox_E\LL(M)^*\subset  (E\ox_{\bbQ_p}L)\ox_E \LL(M)^*$ inside $\calL^{-1}\ox_{E\ox_{\bbQ_p}L}W_{M,\calL}'$. In particular, it fits into a $\GL_2(L)$-equivariant exact sequence 
\begin{align}
    0\to \calL^{-1}\ox_{E\ox_{\bbQ_p}L}\calO[M]\to W_{M,\calL}\to \LL(M)^*\to 0.
\end{align}
Besides, we also have a $\GL_2(L)$-equivariant exact sequence 
\begin{align}
    0\to \calO_\sigma[M]\to \Omega^1_\sigma[M]\to M_{\dR,\sigma}\ox_E \LL(M)^*\to 0
\end{align}
of $E$-vector spaces by taking the $\sigma\in\Sigma$ component of (\ref{eq:dRM}). Here 
\begin{align}
    \calO_\sigma[M]&:=\Hom_{E[D_L^\times]}(\JL(M),E\ox_{\sigma,L}H^0(\calO_{\Dr,\infty}))^{\Gal(\bar L/L)},\\
    \Omega^1_\sigma[M]&:=\Hom_{E[D_L^\times]}(\JL(M),E\ox_{\sigma,L}H^0(\Omega^1_{\Dr,\infty}))^{\Gal(\bar L/L)}.
\end{align}
Let $E_\sigma$ be an $L$-algebra such that $L$ acts via the embedding $\sigma$, and $\calL_\sigma=\calL\ox_{L\ox_{\bbQ_p}E}E_\sigma$ be the $E$-line in $M_{\dR,\sigma}$. Let $\calL_\sigma^{-1}=\calL^{-1}\ox_{L\ox_{\bbQ_p}E}E_\sigma$ be the quotient of $M_{\dR,\sigma}$ such that $\calL_\sigma^{-1}\ox_E\calL_\sigma\isom E$. Define $W_{M,\calL_\sigma}$ to be the inverse image of $\LL(M)^*\subset \calL_\sigma^{-1}\ox_E M_{\dR,\sigma}\ox_E \LL(M)^*$ inside $\calL_\sigma^{-1}\ox \Omega_\sigma^1[M]$. Then one can check that 
\begin{align}
    W_{M,\calL}^*\isom \bigoplus_{\LL(M),\sigma\in\Sigma}W_{M,\calL_\sigma}^*
\end{align}
where $W_{M,\calL_\sigma}^*$ is the strong dual of $W_{M,\calL_\sigma}$, and $W_{M,\calL_\sigma}^*$ is a locally $\sigma$-analytic representation of $\GL_2(L)$ over $E$. In \cite[Proposition 2.10]{CDN20}, they showed that $V_{M,\calL}$-isotypic components of $H^1_{\pro\et}[M]$ are given by $W_{M,\calL}$.

\begin{theorem}\label{thm:VMLmultiplicity}
There is a $\GL_2(L)$-equivariant isomorphism 
\begin{align}
    W_{M,\calL}\isom \Hom_{E[\Gal(\bar L/L)]}(V_{M,\calL},H^1_{\pro\et}[M]).
\end{align}
\end{theorem}
As the locally $\bbQ_p$-analytic representation $W_{M,\calL}^*$ is important, we rename it by $\Pi^{\an}_{\geo}(V_{M,\calL})$. Now, we want to show that $\Pi^{\an}_{\geo}(V_{M,\calL})$ contains exactly the same information of $V_{M,\calL}$. From the discussion in the previous section, as $\WD(M)$ is irreducible, we see the key point is to recover the Hodge filtration. The following theorem gives a direct link between the extension class inside $\Pi^{\an}_{\geo}(V_{M,\calL})^{\sigma\-\lan}$ and the $\sigma$-Hodge filtration of $V_{M,\calL}$.
\begin{theorem}\label{thm:Ext1Dr}
For any $\sigma\in \Sigma$, there exists an isomorphism of $E$-vector spaces
\begin{align}\label{eq:Ext1Dr}
    \WD(M)[-1]\aisom \Ext^1_{\GL_2(L)}(\calO_\sigma[M]^*,\LL(M))
\end{align}
such that via the natural isomorphism 
\begin{align}
    D_{\dR,\sigma}(V_{M,\calL})\isom  \WD(M)[-1]\aisom \Ext^1_{\GL_2(L)}(\calO_\sigma[M]^*,\LL(M))
\end{align}
the $\sigma$-Hodge filtration $\Fil^0 D_{\dR,\sigma}(V_{M,\calL})$ corresponds to the $E$-line generated by $W_{M,\calL_\sigma}^*$. 
\begin{proof}
First we rewrite the above isomorphism on the dual side. From \cite[Th\'eor\`eme 0.4]{CDN20}, we know there is an isomorphism of $E$-vector spaces $M_{\dR,\sigma}\aisom \Hom_{D(G)}(\LL(M)^*,(E\ox_{\sigma,L}H^1_{\dR}(\calM_{\Dr,\infty}))^{\Gal(\bar L/L)})$. Besides, by Theorem \ref{thm:Phiisom} and Remark \ref{rmk:Phiisom2}, the connection map
\begin{align}
    \Hom_{D(G)}(\LL(M)^*,M_{\dR,\sigma}\ox_E \LL(M)^*)\to \Ext^1_{D(G)}(\LL(M)^*,\calO_\sigma[M])
\end{align}
induced by $\Omega^1_\sigma[M]=[\calO_\sigma[M]-M_{\dR,\sigma}\ox_E \LL(M)^*]$ is an isomorphism. (See Remark \ref{rmk:comparetwoExt1map} for more details.) Hence we get an isomorphism of $E$-vector spaces
\begin{align}
    \WD(M)[-1]\aisom M_{\dR,\sigma}\aisom \Hom_{D(G)}(\LL(M)^*,M_{\dR,\sigma}\ox_E \LL(M)^*)\aisom \Ext^1_{D(G)}(\LL(M)^*,\calO_\sigma[M]).
\end{align}
Here we implicitly identified $M_{\dR}$ with $(M[-1])_{\dR}$.

Let $V_{M,\calL}$ be as before, which is a $2$-dimensional $E$-linear de Rham representation of parallel Hodge-Tate weight $0,1$, and $D_{\pst}(V_{M,\calL})\isom M[-1]$. Using the natural identification $D_{\pst}(V_{M,\calL})_\sigma=D_{\dR,\sigma}(V_{M,\calL})$ induced by $B_{\cris}\subset B_{\dR}$ (as $V_{M,\calL}$ is potentially crystalline), we see that the induced map 
\begin{align}
    D_{\dR,\sigma}(V_{M,\calL})\aisom \WD(M)[-1]\aisom \Ext^1_{D(G)}(\LL(M)^*,\calO_\sigma[M])
\end{align}
sends $\calL_\sigma=\Fil^0D_{\dR,\sigma}(V_{M,\calL})$ to the $E$-line generated by the extension class obtained from pulling back $\Omega^1_\sigma[M]$ via $\calL_\sigma\subset D_{\dR,\sigma}(V_{M,\calL})$, which is exactly $W_{M,\calL_\sigma}$ (up to a twist by $E$-line $\calL^{-1}$). This finishes the proof of the theorem.
\end{proof}

\end{theorem}

Here we note the isomorphism in Theorem \ref{thm:Ext1Dr}
\begin{align}
    \WD(M)[-1]\aisom \Ext^1_{\GL_2(L)}(\calO_\sigma[M]^*,\LL(M))
\end{align}
only depends on $M$ and a choice of $\sigma\in\Sigma$.

\begin{remark}\label{rmk:comparetwoExt1map}
Here we breifly discuss the relations between the two isomorphisms between the Weil-Deligne representations and the $\Ext^1$ group, defined in Theorem \ref{thm:Ext1Dr} and in \cite[Theorem 6.3.3]{QS24}. Write $\pi_p=\LL(M)$, $\pi_c=\calO_\sigma[M]^*$ and $\tilde{\pi}=\Omega^1_\sigma[M]^*$ for short. Put $\tilde{\pi}^{\sm}=\pi_0\isom \pi_p^{\oplus 2}$. Then there is a commutative diagram of $E$-vector spaces:
\[\begin{tikzcd}
	{\Hom(\pi_0,\pi_p)} \\
	&& {\Ext^1(\pi_c,\pi_p)} \\
	{\Hom(\pi_p,\pi_0)=\Hom(\pi_p,\tilde{\pi})}
	\arrow["{-\cup[\tilde\pi]}", from=1-1, to=2-3]
	\arrow["{\alpha\mapsto[\ker\alpha\to \pi_0]}"', shift right=3, from=1-1, to=3-1]
	\arrow["{\beta\mapsto[\pi_0\to \coker\beta]}"', from=3-1, to=1-1]
	\arrow["{\gamma\mapsto\coker\gamma}"', from=3-1, to=2-3]
\end{tikzcd}\]
where the equality follows $\Hom(\pi_p,\pi_0)=\Hom(\pi_p,\tilde{\pi})$ from $\Hom(\pi_p,\pi_c)=0$ (by Proposition \ref{prop:Hom}). One directly check that the diagram commutes, and the composition of the two vertical maps are the identity map. We used the above map in Theorem \ref{thm:Ext1Dr} and the below map in \cite[Theorem 6.3.3]{QS24}. If we denote $\Hom(\pi_p,\tilde{\pi})$ by $M$ with a basis $e_1,e_2$, then the map $\beta\mapsto [\pi_0\to \coker\beta]$ is given by $\calL\mapsto \calL^{\perp}=(M/\calL)^*$, and is explicitly given by the $e_1\mapsto e_2^*,e_2\mapsto -e_1^*$. See \cite[Remark 1.2]{Din22} for some related discussions.
\end{remark}

\begin{remark}
The extension group we used in (\ref{eq:Ext1Dr}) is defined as follows. In a recent work \cite[\S 6.2]{QS24} on describing locally analytic vectors in the completed cohomology of unitary Shimura curves (generalizing work of Lue Pan \cite{Pan22,PanII}), we realize the dual of $\calO_\sigma[M]$ inside the completed cohomology of unitary Shimura curves. This establishes the admissibility of $\calO_\sigma[M]$. From this we know that $\calO_\sigma[M]^*$ is an admissible locally $\sigma$-analytic representation of $\GL_2(L)$ over $E$. Therefore, one just put 
\begin{align}
    \Ext^1_{\GL_2(L)}(\calO_\sigma[M]^*,\LL(M)):=\Ext^1_{D_\sigma(\GL_2(L),E)}(\LL(M)^*,\calO_\sigma[M]).
\end{align}
By \cite[Lemme 2.1.1]{Bre19}, the extension group $\Ext^1_{\GL_2(L)}(\calO_\sigma[M]^*,\LL(M))$ is the $E$-vector space of locally $\sigma$-analytic extensions of $\calO_\sigma[M]^*$ by $\LL(M)$. Moreover, from the proof of Corollary \ref{cor:Ext1onlyonDr} we see that the result remains the same if we calculate this extension group in the category of locally $\bbQ_p$-analytic representations.
\end{remark}

Finally, by Proposition \ref{prop:isompext} and Theorem \ref{thm:Ext1Dr}, we get the following result.
\begin{corollary}\label{cor:MainDr}
Let $\rho_p$ be a $2$-dimensional $E$-linear representation of $\Gal(\bar L/L)$, which is de Rham with parallel Hodge-Tate weight $0, 1$, and $D_{\pst}(\rho_p)\isom M[-1]$. Then the $E$-linear isomorphism class of the locally $\bbQ_p$-analytic representation $\Pi^{\an}_{\geo}(\rho_p)$ determines the Hodge filtration of $\rho_p$.
\begin{proof}
Write $\rho_p\isom V_{M,\calL}$ for some $E\ox_{\bbQ_p}L$-line $\calL$ inside $M_{\dR}$. First of all, the smooth vectors inside $\Pi^{\an}_{\geo}(\rho_p)$ determines $M$ by the classical local Langlands correspondence. By Theorem \ref{thm:VMLmultiplicity}, $\Pi^{\an}_{\geo}(\rho_p)$ is an admissible locally $\bbQ_p$-analytic representation of $\GL_2(L)$ over $E$, with its locally $\sigma$-analytic vectors given by a non-split extension $\LL(M)-\calO_\sigma[M]^*$. Here the non-splitness follows from Lemma \ref{lem:nonsplitttt}. By Theorem \ref{thm:Ext1Dr} and Proposition \ref{prop:isompext}, we may recover $\calL_\sigma$ from the $E$-linear isomorphism class of $\Pi^{\an}_{\geo}(\rho_p)^{\sigma\-\lan}$. Finally, as $\sigma$ varies in $\Sigma$, we see the $E$-linear isomorphism class of $\Pi^{\an}_{\geo}(\rho_p)$ determines the Hodge filtration of $\rho_p$.
\end{proof}
\end{corollary}

\bibliographystyle{amsalpha}
\bibliography{DCWC.bib}
\end{document}